\documentclass[10pt,reqno]{amsart}
\usepackage{fixltx2e}                    %% Get the latest LaTeX enhancements
\usepackage[utf8]{inputenc}            %% Customize to your encoding !!!
\usepackage{amsmath}                     %% Imprescindible
\usepackage{amssymb, latexsym, stmaryrd, amsthm, dsfont, amsfonts, amsbsy,amsthm, amsmath, mathrsfs}            %% Imprescindible
\usepackage{mathtools}                   %% for some math typesetting tricks
\usepackage{bm}                          %% for good bold math symbols
\usepackage{enumerate}                   %% enhanced "enumerate" environments
\usepackage{verbatim}                    %% for enhanced "comment" environment
\usepackage{url}                         %% for \url command in bibliography
\usepackage[osf]{mathpazo}  % This provides the palatino font and some nice "old style" numbers.

\usepackage{microtype}  
\usepackage[all, knot]{xy}
        \xyoption{arc} 
        \xyoption{web}                 %% for the best typographic result
\makeatletter                            %% A correction needed for microtype:
\def\MT@register@subst@font{\MT@exp@one@n\MT@in@clist\font@name\MT@font@list
 \ifMT@inlist@\else\xdef\MT@font@list{\MT@font@list\font@name,}\fi}
\makeatother

\usepackage[pdftex,bookmarks,bookmarksnumbered,linktocpage,   %%  customize to your liking
         colorlinks,linkcolor=blue,citecolor=blue]{hyperref}

\newcommand{\bit}{\begin{itemize}}    % but see also \benbullet below
\newcommand{\eit}{\end{itemize}}
\newcommand{\ben}{\begin{enumerate}}
\newcommand{\een}{\end{enumerate}}
\newcommand{\benormal}{\ben[\normalfont 1.]}   % *
\let\enormal\een
\newcommand{\benroman}{\ben[\normalfont (i)]}  % *
\let\eroman\een
\newcommand{\bde}{\begin{description}}
\newcommand{\ede}{\end{description}}

\newcommand{\?}{\ensuremath{\mkern0.4\thinmuskip}}   % very small math space
\let\leq=\leqslant
\let\geq=\geqslant
\let\Box=\square                            % original \Box does not line with the baseline
\let\epsilon=\varepsilon
\let\Lambda\varLambda
\let\Gamma\varGamma
\let\Delta\varDelta
\let\Lambda\varLambda
\let\Omega\varOmega
\let\Theta\varTheta
\let\Xi\varXi
\let\Pi\varPi
\let\Sigma\varSigma
\let\class=\mathsf                              %  classes (of any things)
\let\oper=\mathbb                               %  operators
\bmdefine{\A}{A}                                %  particular algebras
\bmdefine{\2}{2}
\bmdefine{\B}{B}
\bmdefine{\D}{D}
\bmdefine{\M}{M}                                %  the monoid of substitutions
\bmdefine{\LLL}{L}                              %  algebraic language
\bmdefine{\Tm}{Tm}                              %  formula algebra
\bmdefine{\zerou}{[0{,}1]}  
\bmdefine{\T}{T}                                %  particular algebras
\newcommand{\VVV}{\oper{V}}                     %  class operators
\newcommand{\QQQ}{\oper{Q}}
\newcommand{\GGG}{\oper{G}}
\newcommand{\HHH}{\oper{H}}

\newcommand{\PPP}{\oper{P}}
\newcommand{\PSD}{\oper{P}_{\!\textsc{sd}}^{}}
\newcommand{\PPU}{\oper{P}_{\!\textsc{u}}^{}}

\newcommand{\SSS}{\oper{S}}
\newcommand{\III}{\oper{I}}
\newcommand{\UUU}{\oper{U}}

\newcommand{\Con}{\mathrm{Con}}                            % congruences
\newcommand{\KCon}{\mathrm{Comp}}                            % congruences
\bmdefine{\boldstar}{\mathchoice{\textstyle*}{\textstyle*}{\textstyle*}{\scriptstyle*}}

\bmdefine{\btau}{\tau}                                  %  transformer tau
\bmdefine{\brho}{\rho}                                  %  transformer rho

\newcommand{\semeq}{\mathrel{=\joinrel\mathrel\vert\mkern0mu\mathrel\vert\joinrel=}}  % semantical equivalence
\bmdefine{\leibniz}{\Omega}        %  Leibniz operator

\bmdefine{\frege}{\Lambda}         %  Frege operator

\makeatletter
\newcommand{\tarskidsp}{\mathord%
   {\m@th\raisebox{0pt}[0pt][0pt]{$\stackrel%
   {\raisebox{-2.7pt}[0ex][0pt]{$\displaystyle \,\?\thicksim$}}%
   {\displaystyle\leibniz}$}}}
\newcommand{\tarskitxt}{\mathord%
   {\m@th\raisebox{0pt}[0pt][0pt]{$\stackrel%
   {\raisebox{-2.7pt}[0ex][0pt]{$\,\?\thicksim$}}{\displaystyle\leibniz}$}}}
\newcommand{\tarskiscr}{\mathord%
   {{\m@th\raisebox{0pt}[0pt][0pt]{$\stackrel%
   {\raisebox{-2.4pt}[0ex][0pt]{$\scriptstyle \,\?\thicksim$}}%
   {\scriptstyle\leibniz}$}}}}
\newcommand{\tarskiscrscr}{\mathord%
   {{\m@th\raisebox{0pt}[0pt][0pt]{$\stackrel%
   {\raisebox{-2pt}[0ex][0pt]{$\scriptscriptstyle \,\?\thicksim$}}%
   {\scriptscriptstyle\leibniz}$}}}}
\newcommand{\tarski}{\@ifnextchar ^ %
   {\mathchoice{\tarskidsp\kern-.07em}{\tarskitxt\kern-.07em}%
   {\tarskiscr\kern-.07em}{\tarskiscrscr\kern-.07em}}%
   {\mathchoice{\tarskidsp}{\tarskitxt}{\tarskiscr}{\tarskiscrscr}}}
\makeatother

\theoremstyle{theorem}
\newtheorem{Theorem}{Theorem}[section]
\newtheorem{Lemma}[Theorem]{Lemma}
\newtheorem{Corollary}[Theorem]{Corollary}

\newtheorem{fact}{\textbf{Fact}}[Theorem]

\theoremstyle{definition}
\newtheorem{Definition}[Theorem]{Definition}
\newtheorem{exa}[Theorem]{Example}

\theoremstyle{remark}

\newtheorem{Remark}[Theorem]{Remark}

\newcommand{\C}{\boldsymbol{C}} 

\allowdisplaybreaks[1]

\begin{document}
\title[A logical and algebraic characterization of adjunctions]{A logical and algebraic characterization of adjunctions between generalized quasi-varieties}

\author{Tommaso Moraschini}
\email{moraschini@cs.cas.cz}
\address{Institute of Computer Science of the Czech Academy of Science, Prague, Czech Republic}
\date{\today}
\keywords{Adjunction, Adjoint functor, Category Theory, Universal Algebra, Category equivalence, Matrix power, Contextual translation, Locally presentable category\vspace{1mm}}
\maketitle

\begin{abstract}
We present a logical and algebraic description of right adjoint functors between generalized quasi-varieties, inspired by the work of McKenzie on category equivalence. This result is achieved by developing a correspondence between the concept of adjunction and a new notion of translation between relative equational consequences.
\end{abstract}

The aim of the paper is to describe a logical and algebraic characterization of adjunctions between generalized quasi-varieties.\footnote{An extended preprint containing detailed proofs is available on the author's webpage.} This characterization is achieved by developing a correspondence between the concept of adjunction and a new notion of translation, called \textit{contextual translation}, between equational consequences relative to classes of algebras. More precisely, given two generalized quasi-varieties $\class{K}$ and $\class{K}'$, every contextual translation of the equational consequence relative to $\class{K}$ into the one relative to $\class{K}'$ corresponds to a right adjoint functor from $\class{K}'$ to $\class{K}$ and vice-versa (Theorems \ref{Thm : Tr > RA} and \ref{Lemma : RA > TR}). In a slogan, contextual translations between relative equational consequences are the \textit{duals} of right adjoint functors. Examples of this correspondence abound in the literature, e.g., G\"odel's translation of intuitionistic logic into the modal system $\mathcal{S}4$ corresponds to the functor that extracts the Heyting algebra of open elements from an interior algebra (Examples \ref{Exa : God Tr} and \ref{Exa : Open Elements}), and Kolmogorov's translation of classical logic into intuitionistic logic corresponds to the functor that extracts the Boolean algebra of regular elements out of a Heyting algebra (Examples \ref{Exa : Kol Tr} and \ref{Exa : Open Elements}).

The algebraic aspect of our characterization of adjunctions is inspired by the work of McKenzie on category equivalences \cite{Mc96}. Roughly speaking, McKenzie discovered a combinatorial description of category equivalence between prevarieties of algebras (here Theorem \ref{Thm : Decomposition Equivalences}). In particular, he showed that if two prevarieties $\class{K}$ and $\class{K}'$ are categorically equivalent, then we can transform $\class{K}$ into $\class{K}'$ by applying two kinds of deformations to $\class{K}$. The first of these deformations is the \textit{matrix power} construction. The matrix power with exponent $n \in \omega$ of an algebra $\A$ is a new algebra $\A^{[n]}$ with universe $A^{n}$ and whose basic $m$-ary operations are all $n$-sequences of $(m \times n)$-ary term functions of $\A$, which are applied component-wise. The other basic deformation is defined by means of \textit{idempotent and invertible} terms, see Example \ref{Ex : Cat Eq}. This algebraic approach to category equivalences has been reformulated in categorical terms in \cite{HEPorst01AU,HEPorst01} and has an antecedent in \cite{Dukram88}.

Building on McKenzie's work and on the theory of \textit{locally presentable} categories \cite{AdRo94}, we show that every right adjoint functor between generalized quasi-varieties can be decomposed into a combination of two deformations that generalize the ones devised in the special case of category equivalence. These deformations are matrix powers with (possibly) \textit{infinite} exponent and the following construction. Given an algebra $\A$, we say that a set of equations $\theta$ in a single variable is \textit{compatible} with a sublanguage $\mathscr{L}$ of the language of $\A$ if the set of solutions of $\theta$ in $\A$ is closed under the restriction of the operations in $\mathscr{L}$. In this case we let $\theta_{\mathscr{L}}(\A)$ be the algebra obtained by equipping the set of solutions of $\theta$ in $\A$ with the restriction of the operations in $\mathscr{L}$. The main result of the paper shows that every right adjoint functor between generalized quasi-varieties is, up to a natural isomorphism, a composition of the matrix power construction and of the $\theta_{\mathscr{L}}$ construction  (Theorem \ref{Thm:Main}). Moreover, every functor obtained as a composition of these deformations is indeed a right adjoint. This result can be seen as a purely algebraic formulation of the classical description of adjunctions in categories with a free object, which can be traced back at least to \cite{Fr66} (Remark \ref{Rem:Freyd}).

\section{Algebraic preliminaries}

For information on standard notions of universal algebra we refer the reader to \cite{Be11g,BuSa00,McMcTa87}. Given an algebraic language $\mathscr{L}$ and a set $X$, we denote the set of terms over $\mathscr{L}$ built up with the variables in $X$ by $Tm(\mathscr{L}, X)$, and the corresponding absolutely free algebra by $\Tm(\mathscr{L}, X)$. We also denote the set of equations built up from $X$ by $Eq(\mathscr{L}, X)$. Formally speaking, equations are pairs of terms, i.e., $Eq(\mathscr{L}, X) \coloneqq Tm(\mathscr{L}, X) \times Tm(\mathscr{L}, X)$. When the language $\mathscr{L}$ is clear from the context, we simply write $Tm(X)$, $Eq(X)$ and $\Tm(X)$. Sometimes we write $Tm(\mathscr{L}, \kappa)$ to stress the cardinality $\kappa$ of the set of variables. The same convention applies to equations and term algebras. Given two cardinals $\kappa$ and $\lambda$, we denote their Cartesian product by $\kappa \times \lambda$. We denote the set of natural numbers by $\omega$.

We denote  the class operators of isomorphism, homomorphic images, subalgebras, direct products, (isomorphic copies of) subdirect products and ultraproducts respectively by $\III$, $\HHH$, $\SSS$, $\PPP$, $\PSD$ and $\PPU$. We assume that product-style class operators admit empty set of indexes and give a trivial algebra as a result. We denote algebras by bold capital letters $\A$, $\B$, $\C$, etc.\ (with universes $A$, $B$, $C$, etc.). Given a class of algebras $\class{K}$, we denote its language by $\mathscr{L}_{\class{K}}$.

Given an algebraic language $\mathscr{L}$, a \textit{generalized quasi-equation} $\Phi$ is an expression
\[
\Phi = (\bigwedge_{i \in I}\alpha_{i} \thickapprox \beta_{i}) \to \varphi \thickapprox \psi
\]
where $I$ is a possibly infinite \textit{set} and $\alpha_{i}\thickapprox \beta_{i}$ and $\varphi \thickapprox \psi$ are equations. A \textit{quasi-equation} is a generalized quasi-equation in which the set $I$ is finite. Given an algebra $\A$, we say that a generalized quasi-equation $\Phi$ holds in $\A$, in symbols $\A \vDash \Phi$, if for every assignment $\vec{a} \in A$ we have that
\[
\text{if }\alpha_{i}^{\A}(\vec{a}) = \beta^{\A}_{i}(\vec{a}) \text{ for all }i \in I\text{, then }\varphi^{\A}(\vec{a}) = \psi^{\A}(\vec{a}).
\]
A \textit{prevariety}  is a class of algebras axiomatized by (a class of) arbitrary generalized quasi-equations or, equivalently, a class closed under  $\III, \SSS$ and $\PPP$. A \textit{generalized quasi-variety} is a class of algebras axiomatized by (a set of) generalized quasi-equations whose number of variables is bounded by some infinite cardinal.\footnote{It is worth remarking that both the existence and the non-existence of a prevariety that is not a generalized quasi-variety are consistent (relative to large cardinals)  with von Neumann-Bernays-G\"odel class theory (NGB) with the Axiom of Choice. In fact in NBG the assumption that every prevariety is a generalized quasi-variety is equivalent to the \textit{Vop\v{e}nka Principle}, which states that every class of pairwise non-embeddable models of a first-order theory is a set \cite{Ad90AU} (see also \cite[Proposition 2.3.18]{Go98a}).} These can be equivalently characterized \cite{BlJo06} as the classes of algebras closed under $\III, \SSS, \PPP$ and $\UUU_{\kappa}$ (for some infinite cardinal $\kappa$), where for every class of algebras $\class{K}$,
\[
\UUU_{\kappa}(\class{K}) \coloneqq \{ \A : \B \in \class{K}\text{ for every $\kappa$-generated subalgebra }\B \leq \A \}.
\]
It is well known that a \textit{quasi-variety} is a class of algebras axiomatized by quasi-equations or, equivalently, a class closed under $\III$, $\SSS$, $\PPP$ and $\PPU$. A \textit{variety} is a class of algebras axiomatized by equations or, equivalently, closed under $\HHH$, $\SSS$ and $\PPP$. Given a class of algebras $\class{K}$, we denote by $\GGG\QQQ_{\kappa}(\class{K})$ the models of the generalized quasi-equations in $\kappa$-many variables that hold in $\class{K}$ and respectively by \label{DefOfQQQ}$\QQQ(\class{K})$ and $\VVV(\class{K})$ the quasi-variety and the variety generated by $\class{K}$. It is well known that
\[
\GGG\QQQ_{\kappa}(\class{K}) = \UUU_{\kappa}\III\SSS\PPP(\class{K}) \quad\QQQ(\class{K}) = \III\SSS\PPP\PPU(\class{K}) \quad \VVV(\class{K}) = \HHH\SSS\PPP(\class{K}).
\]

Given a class of algebras $\class{K}$ and a set $X$, we denote by $\Tm_{\class{K}}(X)$ the free algebra in $\class{K}$ with free generators $X$. In general the free algebra $\Tm_{\class{K}}(X)$ is constructed as a quotient of the term algebra $\Tm(X)$ and its elements are congruence classes of terms  equivalent in $\class{K}$. Sometimes we identify the universe of $\Tm_{\class{K}}(X)$ with a set of its representatives, i.e., with a set of terms in variables $X$. It is well known that prevarieties contain free algebras with arbitrary large sets of free generators.

Given a class of algebras $\class{K}$ and an algebra $\A$, we say that a congruence $\theta$ of $\A$ is a $\class{K}$-\textit{congruence} if $\A / \theta \in \class{K}$, and 
denote the collection of $\class{K}$-congruences by $\Con_{\class{K}}\A$. In particular, we will denote by $\pi_{\theta} \colon \A \to \A / \theta$ the canonical surjection. If $\class{K}$ is a prevariety, then $\Con_{\class{K}}\A$ forms a closure system when ordered under the inclusion relation.  
We  denote by $\textup{Cg}_{\class{K}}^{\A}$ the closure 
operator of generation  of $\class{K}$-congruences.

Given a class of algebras $\class{K}$ and $\Phi \cup \{ \epsilon \thickapprox \delta \} \subseteq Eq(X)$, we define
\begin{align*}
\Phi \vDash_{\class{K}} \epsilon \thickapprox \delta \Longleftrightarrow &\text{ for every }\A \in \class{K} \text{ and every }h \colon \Tm(X) \to \A\\
&\text{ if }h \varphi  = h \psi \text{ for every }\varphi \thickapprox \psi \in \Phi \text{, then }h \epsilon  = h \delta.
\end{align*}
The relation $\vDash_{\class{K}}$ is called the \textit{equational consequence relative to} $\class{K}$. The function $C_{\class{K}} \colon \mathcal{P}(Eq(X)) \to \mathcal{P}(Eq(X))$ defined by the rule
\[
C_{\class{K}}(\Phi) \coloneqq \{ \epsilon \thickapprox \delta : \Phi \vDash_{\class{K}} \epsilon \thickapprox \delta \}\text{, for every }\Phi \subseteq Eq(X)
\]
is a closure operator over $Eq(X)$. If $\class{K}$ is a prevariety, then the set of fixed points of $C_{\class{K}} \colon \mathcal{P}(Eq(X)) \to \mathcal{P}(Eq(X))$ coincides with $\Con_{\class{K}}\Tm(X)$. Now let $\class{K}$ be a quasi-variety and $\A$ an arbitrary algebra. The lattice $\Con_{\class{K}}\A$ is algebraic and its compact elements $\KCon_{\class{K}}\A$ are the finitely 
generated $\class{K}$-congruences. In particular, the closure operator $\textup{Cg}_{\class{K}}^{\A}$ is finitary. An algebra $\A \in \class{K}$ is $\class{K}$-\textit{finitely presentable} if there is some $n \in \omega$ and $\theta \in \KCon_{\class{K}}\Tm_{\class{K}}(n)$ such that $\A$ is isomorphic to $\Tm_{\class{K}}(n) / \theta$. 

For standard information on category theory we refer the reader to \cite{AdHeSt06,Aw06,La98c}, while for categorical universal algebra see \cite{AdRo94,AdRoVi11}. Let $\kappa$ be a regular cardinal and let $\class{K}$ be a locally small category. An object $\A$ in $\class{K}$ is $\kappa$\textit{-presentable} if the functor $\hom(\A, \cdot)$ preserves $\kappa$-directed colimits. In generalized quasi-varieties the $\kappa$-presentable objects can be described as follows:

\begin{Lemma}\label{Lem : What are presetable}
Let $\kappa$ be a regular cardinal and $\class{K}$ be a generalized quasi-variety axiomatized by generalized quasi-equations in less than $\kappa$ variables. An algebra $\A \in \class{K}$ is $\kappa$-presentable in the categorical sense if and only if it is (isomorphic to) a quotient of $\Tm_{\class{K}}(\lambda)$ under a $\mu$-generated $\class{K}$-congruence for some $\lambda, \mu < \kappa$.
\end{Lemma}

Let $\kappa$ be a regular cardinal and $\class{K}$ be a locally small category. $\class{K}$ is \textit{locally $\kappa$-presentable} \cite{AdRo94} if it is cocomplete, and has a \textit{set} $J$ of $\kappa$-presentable objects such that every object in $\class{K}$ is a $\kappa$-directed colimit of objects in $J$. Moreover, $\class{K}$ is \textit{locally presentable} if it is locally $\kappa$-presentable for some regular cardinal $\kappa$. Generalized quasi-varieties, equipped with homomorphisms, can be seen as categories.\footnote{If the language of a generalized quasi-variety $\mathsf{K}$ contains no constant symbols, then (when seen as a category) $\mathsf{K}$ is assumed to contain the \textit{empty algebra} as an object.} 

\begin{Lemma}\label{Lem:LocallyPresentable}
Generalized quasi-varieties are locally presentable categories.
\end{Lemma}

Ad\'amek and Rosick\'{y} proved in \cite[Theorem 1.66]{AdRo94} the following characterization of right adjoint functors between locally presentable categories. By Lemma \ref{Lem:LocallyPresentable} it applies to generalized quasi-varieties as well.
\begin{Theorem}[Ad\'amek and Rosick\'{y}]\label{Thm:Adamek}
A functor between locally presentable categories is right adjoint if and only if it preserves limits and $\kappa$-directed colimits for some regular cardinal $\kappa$.
\end{Theorem}

Given two prevarieties $\class{X}$ and $\class{Y}$, the functors $\mathcal{F} \colon \class{X} \longleftrightarrow \class{Y} \colon \mathcal{G}$, where $\mathcal{F}$ sends everything to the initial object and $\mathcal{G}$ sends every object to the terminal object, always form an adjunction $\mathcal{F} \dashv \mathcal{G}$. We call the adjunctions of this kind \textit{trivial}. In particular, we say that a left (resp.\ right) adjoint functor between prevarieties is \textit{trivial} if it sends everything to the initial (resp.\ terminal) object.

\section{The two basic deformations}\label{Sec:Constructions}

In this section we describe two general methods to deform a given generalized quasi-variety, obtaining a new generalized quasi-variety that is related to the first one by an adjunction. The first deformation that we consider is just an infinite version of the usual \textit{finite matrix power} construction. Let $\class{X}$ be a class of similar algebras and $\kappa$ be a cardinal. Then observe that every term $\varphi \in Tm(\kappa)$ induces a map $\varphi \colon A^{\kappa} \to A$ for every $\A \in \class{X}$.

\begin{Definition}\label{Def:Chap5:Language}
Let $\kappa > 0$ be a cardinal and $\class{X}$ a class of similar algebras. Then $\mathscr{L}_{\class{X}}^{\kappa}$ is the algebraic language whose $n$-ary operations (for every $n \in \omega$)  are all $\kappa$-sequences $\langle  t_{i} : i < \kappa \rangle$ of terms $t_{i}$ of the language of $\class{X}$ built up with variables
\[
\{ x^{j}_{m} : 1 \leq m \leq  n \text{ and }j < \kappa \}.
\]
Observe that each $t_i$ has a finite number of variables, possibly none, of each sequence $\vec{x}_{m} \coloneqq \langle x^{j}_{m} : j < \kappa \rangle$ with $1 \leq m \leq n$. We will write $t_{i} = t_{i}(\vec{x}_{1}, \dots, \vec{x}_{n})$ to denote this fact.
\end{Definition}

\begin{exa}\label{Exa:MatrixLanguage}
Consider the variety of bounded distributive lattices $\class{DL}_{01}$. Examples of basic operations of $\mathscr{L}^{2}_{\class{DL}_{01}}$ are:
\begin{align*}
\langle x^{1}, x^{2}\rangle \sqcap \langle y^{1}, y^{2}\rangle &\coloneqq  \langle x^{1} \land y^{1}, x^{2} \lor y^{2}\rangle\\
\langle x^{1}, x^{2}\rangle \sqcup \langle y^{1}, y^{2}\rangle &\coloneqq \langle x^{1} \lor y^{1}, x^{2} \land y^{2}\rangle\\
\lnot \langle x^{1}, x^{2}\rangle &\coloneqq \langle x^{2}, x^{1}\rangle\\
1 &\coloneqq \langle 1, 0 \rangle\\
0 &\coloneqq  \langle 0, 1 \rangle. 
\end{align*}
\end{exa}

\begin{Definition}\label{Def:Chap5:MatrixPower}
Consider an algebra $\A \in \class{X}$ and a cardinal $\kappa > 0$. We let $\A^{[\kappa]}$ be the algebra of type $\mathscr{L}_{\class{X}}^{\kappa}$ with universe $A^{\kappa}$ where a $n$-ary operation $\langle t_{i} : i < \kappa \rangle$ is interpreted as 
\[
\langle t_{i} : i < \kappa \rangle( a_{1}, \dots, a_{n}) = \langle t_{i}^{\A}(a_{1} / \vec{x}_{1}, \dots, a_{n}/\vec{x}_{n}): i < \kappa \rangle
\]
for every $a_{1}, \dots, a_{n} \in A^{\kappa}$ (the notation $a_{m} / \vec{x}_{m}$ means that we are assigning the tuple $a_{m}$ of elements of $A$ to the tuple of variables $\vec{x}_{m}$). In other words $\langle t_{i} : i < \kappa \rangle( a_{1}, \dots, a_{n})$ is the $\kappa$-sequence of elements of $A$ defined as follows. Consider $i < \kappa$. Observe that only a finite number of variables occurs in $t_{i}$, say
\[
t_{i} = t_{i}(x_{1}^{\alpha_{1}^{1}}, \dots, x_{1}^{\alpha_{m_{1}}^{1}}, \dots, x_{n}^{\alpha_{1}^{n}}, \dots, x_{n}^{\alpha_{m_{n}}^{n}}),
\]
where $\alpha_{1}^{1}, \dots, \alpha_{m_{1}}^{1}, \dots, \alpha_{1}^{n}, \dots, \alpha_{m_{n}}^{n} < \kappa$. Then the $i$-th component of the sequence $\langle t_{i} : i < \kappa \rangle( a_{1}, \dots, a_{n})$ is
\[
t_{i}^{\A}(a_{1}(\alpha_{1}^{1}), \dots, a_{1}(\alpha_{m_{1}}^{1}), \dots, a_{n}(\alpha_{1}^{n}), \dots, a_{1}(\alpha_{m_{n}}^{n})).
\]
If $\class{X}$ is a class of similar algebras, we set
\[
\class{X}^{[\kappa]} \coloneqq \III \{ \A^{[\kappa]} : \A \in \class{X} \}
\]
and call it the $\kappa$-th \textit{matrix power} of $\class{X}$.
\end{Definition}

Now, let $[\kappa]$ be the map defined as follows:
\begin{align*}
\A &\longmapsto \A^{[\kappa]}\\
f \colon \A \to \B &\longmapsto f^{[\kappa]} \colon \A^{[\kappa]} \to \B^{[\kappa]}
\end{align*}
where $f^{[\kappa]} \langle a_{i} : i < \kappa \rangle \coloneqq \langle f(a_{i}) : i < \kappa \rangle$, for every $\A, \B \in \class{X}$ and every homomorphism $f$. It is easy to check that the map $f^{[\kappa]} \colon \A^{[\kappa]} \to \B^{[\kappa]}$ is indeed a homomorphism.

\begin{exa}\label{Exa:HowLanguageBelongs}
In Example \ref{Exa:MatrixLanguage} we highlighted some operations of $\mathscr{L}_{\class{DL}_{01}}^{2}$. Let us explain how are they interpreted in the matrix power construction. Consider $\A \in \class{DL}_{01}$. The universe of $\A^{[2]}$ is just the Cartesian product $A \times A$. We have that:
\begin{align*}
\langle a, b \rangle \sqcap \langle c, d \rangle &= \langle a \land c, b \lor d \rangle\\
\langle a, b \rangle \sqcup \langle c, d \rangle &= \langle a \lor c, b \land d \rangle\\
\lnot \langle a, b \rangle & = \langle b, a \rangle\\
1 &= \langle 1^{\A}, 0^{\A}\rangle\\
0 &= \langle 0^{\A}, 1^{\A}\rangle
\end{align*}
for every $\langle a, b\rangle, \langle c, d\rangle \in A \times A$. Examples of matrix powers with infinite exponent are technically, but not conceptually, more involved (see Example  \ref{Ex:Chap6:Hom-Functor}).
\qed
\end{exa}

\begin{Theorem}\label{Thm:McKenzie}
Let $\class{X}$ be a generalized quasi-variety and $\kappa > 0$ a cardinal. If $\class{Y}$ is a generalized quasi-variety such that $\class{X}^{[\kappa]} \subseteq \class{Y}$, then $[ \kappa ] \colon \class{X} \to \class{Y}$ is a right adjoint functor.
\end{Theorem}

\begin{proof}
It is not difficult to see that the map $[\kappa]$ is a functor that preserves direct products and equalizers. Since all limits can be constructed as combination of products and equalizers, we conclude that $[\kappa]$ preserves limits. Moreover, $[\kappa]$ preserves $\lambda$-directed colimits for every regular cardinal $\lambda$, larger than the number of variables occurring in the generalized quasi-equations axiomatizing $\class{X}$ and $\class{Y}$. With an application of Theorem \ref{Thm:Adamek} we are done.
\end{proof}

\begin{exa}[\textsf{Finite Exponent}]\label{Exa : Matrix}
It is not difficult to see that if $\class{X}$ is a class of similar algebras and $\kappa > 0$, then the functor $[\kappa] \colon \class{X} \to \class{X}^{[\kappa]}$ is a category equivalence (see for example \cite[Theorem 2.3.(i)]{Mc96} where this is stated under the assumption that $\kappa$ is  finite). Moreover, when $\kappa$ is finite, it happens that if $\class{X}$ is a prevariety (or a generalized quasi-variety, a quasi-variety, a variety), then so is $\class{X}^{[\kappa]}$.
\qed
\end{exa}

In order to describe the second kind of deformation, we need the following:

\begin{Definition}
Let $\class{X}$ be a class of similar algebras and $\mathscr{L} \subseteq \mathscr{L}_{\?\class{X}}$. A set of equations $\theta \subseteq Eq( \mathscr{L}_{\?\class{X}}, 1)$ is \textit{compatible} with $\mathscr{L}$ in $\class{X}$ if for every $n$-ary operation $\varphi \in \mathscr{L}$ we have that
\[
\theta(x_{1}) \cup \dots \cup \theta(x_{n}) \vDash_{\class{X}} \theta(\varphi (x_{1}, \dots, x_{n})).
\]
\end{Definition}
\noindent In other words $\theta$ is compatible with $\mathscr{L}$  in $\class{X}$ when the solution sets of $\theta$ in $\class{X}$ are closed under the interpretation of the operations and constants in $\mathscr{L}$.

Now we will explain how is it possible to build a functor out of a set of equations $\theta$ compatible with $\mathscr{L} \subseteq \mathscr{L}_{\?\class{X}}$. For every $\A \in \class{X}$, we let $\theta_{\mathscr{L}}(\A)$ be the algebra of type $\mathscr{L}$ whose universe is
\[
\theta_{\mathscr{L}}(A) \coloneqq \{ a \in A : \A \vDash \theta(a) \}
\]
equipped with the restriction of the operations in $\mathscr{L}$. We know that $\theta_{\mathscr{L}}(\A)$ is well-defined, since its universe is closed under the interpretation of the operations in $\mathscr{L}$ and contains the interpretation of the constants in $\mathscr{L}$. Observe that by definition of compatibility $\theta_{\mathscr{L}}(\A)$ can be empty if and only if $\mathscr{L}$ contains no constant symbol.

Given  a homomorphism $f \colon \A \to \B$ in $\class{X}$, we denote its restriction to $\theta_{\mathscr{L}}(A)$ by
\[
\theta_{\mathscr{L}}(f) \colon  \theta_{\mathscr{L}}(\A) \to \theta_{\mathscr{L}}(\B).
\]
It is easy to see that $\theta_{\mathscr{L}}(f)$ is a well-defined homomorphism. Now, consider the following class of algebras:
\[
\theta_{\mathscr{L}}(\class{X}) \coloneqq \III \{ \theta_{\mathscr{L}}(\A) : \A \in \class{X} \}.
\]
Let $\theta_{\mathscr{L}} \colon \class{X} \to \theta_{\mathscr{L}}(\class{X})$ be the map defined by the following rule:
\begin{align*}
\A & \longmapsto \theta_{\mathscr{L}}(\A)\\
f \colon \A \to \B & \longmapsto\theta_{\mathscr{L}}(f) \colon  \theta_{\mathscr{L}}(\A) \to \theta_{\mathscr{L}}(\B).
\end{align*}
It is easy to check that $\theta_{\mathscr{L}}$ is a functor.

\begin{Theorem}\label{Thm:ThetaRight}
Let $\class{X}$ be a generalized quasi-variety and $\theta \subseteq  Eq( \mathscr{L}_{\?\class{X}}, 1)$ a set of equations compatible with $\mathscr{L} \subseteq \mathscr{L}_{\?\class{X}}$. If $\class{Y}$ is a generalized quasi-variety such that $\theta_{\mathscr{L}}(\class{X}) \subseteq \class{Y}$, then $\theta_{\mathscr{L}} \colon \class{X} \to \class{Y}$ is a right adjoint functor.
\end{Theorem}
\begin{proof}
By Theorem \ref{Thm:Adamek} we know that the functor $\theta_{\mathscr{L}}$ is a right adjoint if and only if it preserves limits and $\kappa$-directed colimits for some regular cardinal $\kappa$. It is easy to see that $\theta_{\mathscr{L}}$ preserves direct products and equalizers and, therefore, all limits. Moreover, $\theta_{\mathscr{L}}$ preserves $\kappa$-directed colimits for every regular cardinal $\kappa$, larger than the number of variables occurring in the generalized quasi-equations axiomatizing $\class{X}$.
\end{proof}

A familiar instance of the above construction is the following:
\begin{exa}[\textsf{Subreducts}]
Let $\class{X}$ be a (generalized) quasi-variety and $\mathscr{L} \subseteq \mathscr{L}_{\?\class{X}}$. An $\mathscr{L}$-\textit{subreduct} of  an algebra $\A \in \class{X}$ is a subalgebra of the $\mathscr{L}$-reduct of $\A$. From \cite[Proposition 2.3.19]{Go98a} it is easy to infer that the class $\class{Y}$ of $\mathscr{L}$-subreducts of algebras in $\class{X}$ is a (generalized) quasi-variety. For quasi-varieties this fact was proved by Maltsev \cite{Ma71}. Consider the forgetful functor $\mathcal{U} \colon \class{X} \to \class{Y}$. It is easy to see that $\mathcal{U} = \theta_{\mathscr{L}}$ where $\theta = \emptyset$. From Theorem \ref{Thm:ThetaRight} it follows that $\mathcal{U}$ has a left adjoint.
\qed
\end{exa}

In the next examples we illustrate how the two deformations introduced so far can be combined to describe right adjoint functors.

\begin{exa}[\textsf{Kleene Algebras}]\label{Exa : KA}
A \textit{Kleene algebra} $\A = \langle A, \sqcap, \sqcup, \lnot, 0, 1\rangle$ is a De Morgan algebra in which the equation $x \sqcap \lnot x \leq y \sqcup \lnot y$ holds. We denote by $\class{KA}$ the variety of Kleene algebras and by $\class{DL}_{01}$ the variety of bounded distributive lattices. In \cite{Ci86} (but see also \cite{Ka58}) a way of constructing Kleene algebras out of bounded distributive lattices is described. More precisely, given $\A \in \class{DL}_{01}$, the Kleene algebra $\mathcal{G}(\A)$ has universe
\[
G(A) \coloneqq \{ \langle a, b\rangle \in A^{2} : a \land b = 0 \}
\]
and operations defined as in Example \ref{Exa:HowLanguageBelongs}. Moreover, given a homomorphism $f \colon \A \to \B$ in $\class{DL}_{01}$, the map $\mathcal{G}(f) \colon \mathcal{G}(\A) \to \mathcal{G}(\B)$ is defined by replicating $f$ component-wise.  It turns out that $\mathcal{G} \colon\class{DL}_{01} \to \class{KA}$ is a right adjoint functor \cite[Theorem 1.7]{Ci86}.

In order to decompose $\mathcal{G}$ into a combination of our two deformations, consider the sublanguage $\mathscr{L}$ of the language of $\class{DL}^{[2]}_{01}$, defined in Example \ref{Exa:MatrixLanguage}. Consider also the set of equations
\[
\theta \coloneqq \{ \langle x^{1}\land x^{2}, x^{1}\land x^{2}\rangle \thickapprox \langle 0, 0 \rangle \} \subseteq Eq(\mathscr{L}_{\? \class{DL}^{[2]}_{01}}, 1).
\]
It is easy to see that $\theta$ is compatible with $\mathscr{L}$. Moreover, for every $\A \in \class{DL}_{01}$ and $a, b \in A$ we have that
\[
\langle a, b \rangle \in \mathcal{G}(\A) \Longleftrightarrow \langle a, b \rangle \in \theta_{\mathscr{L}}(A).
\]
Hence we conclude that $\theta_{\mathscr{L}}(\A^{[2]}) = \mathcal{G}(\A) \in \class{KA}$ for every $\A\in \class{DL}_{01}^{[2]}$. This implies  that the functor $\mathcal{G}$ coincides with the composition $\theta_{\mathscr{L}} \circ [2]$, where $[2]\colon \class{DL}_{01} \to \class{DL}^{[2]}_{01}$.
\qed
\end{exa}

Before concluding this section, we show that the deformations described until now can be applied to decompose equivalence functors between prevarieties. This will make the connection with McKenzie's work \cite{Mc96} explicit. To this end, let us recall the definition of a special version of the $\theta_{\mathscr{L}}$ construction.

\begin{exa}[\textsf{Idempotent and Invertible Terms}]\label{Ex : Cat Eq}
Suppose that $\class{X}$ is a prevariety and $\sigma(x)$ a unary term. We say that $\sigma(x)$ is \textit{idempotent} if $\class{X}\vDash \sigma \sigma(x) \thickapprox \sigma(x)$ and that $\sigma(x)$ is \textit{invertible} if there are an $n$-ary term $t$ and unary terms $t_{1}, \dots, t_{n}$ such that
\[
\class{X} \vDash t( \sigma t_{1}(x), \dots, \sigma t_{n}(x)) \thickapprox x.
\]
Given a unary and idempotent term $\sigma(x)$ of $\class{X}$, we define 
\[
\mathscr{L} \coloneqq \{ \sigma t : t \text{ is a basic symbol of }\class{X}^{[1]}\}
\] 
and $\theta \coloneqq \{ x \thickapprox \sigma(x) \}$. Moreover, we define
\[
\class{X}(\sigma) \coloneqq \theta_{\mathscr{L}}(\class{X}^{[1]}).
\]
McKenzie proved that the functor $\sigma \colon \class{X} \to \class{X}(\sigma)$ defined as the composition $\theta_{\mathscr{L}} \circ [1]$ is a category equivalence \cite[Theorem 2.2.(ii)]{Mc96}. Moreover, if $\class{X}$ is a prevariety (or a generalized quasi-variety, a quasi-variety, a variety), then so is $\class{X}(\sigma)$. Following the literature, we will write $\A(\sigma)$ instead of $\sigma(\A)$ for every $\A \in \class{X}$. 
\qed
\end{exa}

To introduce McKenzie's characterization of category equivalence, we restrict to prevarieties without constant symbols. It should be kept in mind that this restriction is somehow immaterial, since, given a prevariety $\class{K}$, we can always replace the constant symbols of $\class{K}$ by constant unary operations obtaining a new prevariety $\class{K}'$ whose only difference with $\class{K}$ is the presence of the empty algebra.

We need to recall some basic concepts \cite[Definitions 4.76 and 4.77]{Be11g}:

\begin{Definition}\label{Def:Term-Eqq}
Let $\class{X}$ and $\class{Y}$ be prevarieties without constant symbols.
An \textit{interpretation} of $\class{X}$ in $\class{Y}$ is a map $\btau \colon \mathscr{L}_{\class{X}} \to Tm(\mathscr{L}_{\class{Y}}, \omega)$ such that:
\benormal
\item $\btau$ sends $n$-ary basic symbols to at most $n$-ary terms for every $n \geq 1$.
%\item $\btau$ sends constants symbols to at most unary terms.
\item $\A^{\btau} \coloneqq \langle A, \{ \btau(\lambda) : \lambda \in \mathscr{L}_{\class{X}} \}\rangle \in \class{X}$ for every $\A \in \class{Y}$.
\enormal
\end{Definition}

\begin{Definition}
Two prevarieties $\class{X}$ and $\class{Y}$ without constant symbols are \textit{term-equivalent} if there are interpretations $\btau$ and $\brho$ of $\class{X}$ in $\class{Y}$ and of $\class{Y}$ in $\class{X}$ respectively such that for every $\A \in \class{X}$ and $\B \in \class{Y}$,
\[
(\A^{\brho})^{\btau} = \A \text{ and }(\B^{\btau})^{\brho} = \B.
\]
\end{Definition}

When two prevarieties $\class{X}$ and $\class{Y}$ without constant symbols are term-equivalent, the map  that sends $\A \in \class{X}$ to $\A^{\brho} \in \class{Y}$ and that is the identity on arrows is a category equivalence $\mathcal{F}^{\brho} \colon \class{X} \to \class{Y}$. Then we have the following \cite[Theorem 6.1]{Mc96}:
\begin{Theorem}[McKenzie]\label{Thm : Decomposition Equivalences}
If $\mathcal{G} \colon \class{X} \to \class{Y}$ is a category equivalence between prevarieties without constant symbols, then there are a natural number $n > 0$ and a unary idempotent and invertible term $\sigma(x)$ of $\class{X}^{[n]}$ such that
\benormal
\item $\class{Y}$ is term-equivalent to $\class{X}^{[n]}(\sigma)$ under some interpretation $\brho$ of $\class{Y}$ in $\class{X}^{[n]}(\sigma)$.
\item The functors $\mathcal{G}$ and $\mathcal{F}^{\brho}\circ (\sigma \circ [n])$ are naturally isomorphic.
\enormal
\end{Theorem}

\section{From translations to right adjoints}\label{Sec : TR}

As we mentioned, our aim is to develop a correspondence between the \textit{adjunctions} between two generalized quasi-varieties $\class{X}$ and $\class{Y}$ and the \textit{translations} between the equational consequences relative to $\class{X}$ and $\class{Y}$. To simplify the notation, we will assume throughout this section that $\class{X}$ and $\class{Y}$ are two fixed generalized quasi-varieties (possibly in different languages).

\begin{Definition}
Consider a cardinal $\kappa > 0$. A $\kappa$-\textit{translation} $\btau$ of $\mathscr{L}_{\?\class{X}}$ into $\mathscr{L}_{\?\class{Y}}$ is a map from $\mathscr{L}_{\?\class{X}}$ to $\mathscr{L}_{\class{Y}}^{\kappa}$ that preserves the arities of function symbols.
\end{Definition}
\noindent In other words, if a basic symbol $\varphi \in \mathscr{L}_{\?\class{X}}$ is $n$-ary, we have that $\btau(\varphi) = \langle t_{i} : i < \kappa \rangle$ for some terms $t_{i} = t_{i}( \vec{x}_{1}, \dots, \vec{x}_{n})$ of language of $\class{Y}$, where $\vec{x}_{m}= \langle x^{j}_{m} : j < \kappa \rangle$. It is worth  remarking that $\btau$ sends constant symbols to sequences of constant symbols. Thus if $\mathscr{L}_{\?\class{X}}$ contains a constant symbol, then also $\mathscr{L}_{\?\class{Y}}$ must contain one for a translation to exist.

A $\kappa$-translation $\btau$ extends naturally to arbitrary terms. Let us explain briefly how. Given a cardinal $\lambda$, let $Tm(\mathscr{L}_{\?\class{X}}, \lambda)$ be the set of terms of $\class{X}$ written with variables in $\{ x_{j} : j < \lambda \}$ and let $Tm(\mathscr{L}_{\?\class{Y}}, \kappa \times \lambda)$ be the set of terms of $\class{Y}$ written with variables in $\{ x_{j}^{i} : j < \lambda, i < \kappa \}$. We define recursively a map
\[
\btau_{\ast} \colon Tm(\mathscr{L}_{\?\class{X}}, \lambda) \to Tm(\mathscr{L}_{\?\class{Y}}, \kappa \times \lambda)^{\kappa}.
\]
For variables and constants we set
\begin{align*}
\btau_{\ast}(x_{j}) &\coloneqq\langle x^{i}_{j} : i < \kappa \rangle, \text{ for every }j < \lambda\\
\btau_{\ast}(c) & \coloneqq \btau(c).
\end{align*}
For complex terms, let $\psi \in \mathscr{L}_{\? \class{X}}$ be $n$-ary and $\varphi_{1}, \dots, \varphi_{n} \in Tm(\mathscr{L}_{\?\class{X}}, \lambda)$. We have that $\btau(\psi) = \langle t_{i} : i < \kappa \rangle$ where $t_{i} = t_{i}(\vec{x}_{1}, \dots, \vec{x}_{n})$. Keeping this in mind, we set
\[
\btau_{\ast}(\psi(\varphi_{1}, \dots, \varphi_{n}) )(i) \coloneqq t_{i}( \btau_{\ast}(\varphi_{1}) / \vec{x}_{1}, \dots, \btau_{\ast}(\varphi_{n}) / \vec{x}_{n})\text{ for every }i < \kappa.
\]

The map $\btau_{\ast}$ can be lifted to sets of equations yielding a new function
\[
\btau^{\ast} \colon \mathcal{P} ( Eq(\mathscr{L}_{\?\class{X}}, \lambda) ) \to \mathcal{P} ( Eq(\mathscr{L}_{\?\class{Y}}, \kappa \times \lambda) )
\]
defined by the following rule:
\[
\Phi \longmapsto \{ \btau_{\ast}(\epsilon)(i) \thickapprox \btau_{\ast}(\delta)(i) : i < \kappa \text{ and }\epsilon \thickapprox \delta \in \Phi \}.
\]
\begin{Definition}\label{Def : Translation}
Consider a cardinal $\kappa > 0$. A \textit{contextual} $\kappa$-\textit{translation} of $\vDash_{\class{X}}$ into $\vDash_{\class{Y}}$ is a pair $\langle \btau, \Theta \rangle$ where $\btau$ is a $\kappa$-translation of $\mathscr{L}_{\?\class{X}}$ into $\mathscr{L}_{\?\class{Y}}$ and $\Theta(\vec{x}) \subseteq Eq(\mathscr{L}_{\?\class{Y}}, \kappa)$ is a set of equations written with variables among $\{ x^{i}: i < \kappa \}$ that satisfies the following conditions:

\benormal
\item For every cardinal $\lambda$ and equations $\Phi \cup \{ \epsilon \thickapprox \delta \} \subseteq Eq(\mathscr{L}_{\?\class{X}}, \lambda)$ written with variables among $\{ x_{j} : j < \lambda \}$,
\[
\text{if }\Phi \vDash_{\class{X}}  \epsilon \thickapprox \delta\text{, then }\btau^{\ast}(\Phi) \cup \bigcup_{j < \lambda} \Theta (\vec{x}_{j}) \vDash_{\class{Y}} \btau^{\ast}(\epsilon \thickapprox \delta).
\]
\item For every $n$-ary operation $\psi \in \mathscr{L}_{\? \class{X}}$,
\[
\Theta (\vec{x}_{1}) \cup \dots \cup \Theta (\vec{x}_{n}) \vDash_{\class{Y}} \Theta ( \btau_{\ast}\psi(x_{1}, \dots, x_{n})).
\]
\enormal
In 1 and 2 it is intended that $\vec{x}_{j} = \langle x^{i}_{j} : i < \kappa \rangle$. The set $\Theta$ is the \textit{context} of the contextual translation $\langle \btau, \Theta \rangle$.

A contextual $\kappa$-translation $\langle \btau, \Theta \rangle$ of $\vDash_{\class{X}}$ into $\vDash_{\class{Y}}$ is \textit{non-trivial} provided that if there is a (non-empty) sequence $\vec{\varphi} \in Tm(\mathscr{L}_{\?\class{Y}}, 0)^{\kappa}$ of constant symbols such that $\class{Y} \vDash \Theta (\vec{\varphi})$, then there is $i_{0} < \kappa$ and sequences of variables
\[
\vec{x} = \langle x^{i} : i < \kappa \rangle \text{ and }\vec{y} = \langle y^{i} : i < \kappa \rangle
\]
such that
\[
\Theta ( \vec{x}) \cup \Theta ( \vec{y}) \nvDash_{\class{Y}} x^{i_{0}} \thickapprox y^{i_{0}}.
\]
\end{Definition}

Several translations between logics classically considered in the literature provide examples of this general notion of contextual translation. 

\begin{exa}[\textsf{Heyting and Interior Algebras}]\label{Exa : God Tr}
As shown by G\"odel in \cite{Go33b} (see also \cite{DuLe59,MakRyb74,McKT48}), it is possible to interpret the intuitionistic propositional calculus $\mathcal{IPC}$ into the consequence relation associated with the global modal system $\mathcal{S}4$ \cite{Kr99,MK07c}. Since these two logics are algebraizable \cite{BP89} with equivalent algebraic semantics the variety of Heyting algebras $\class{HA}$ and of interior algebras $\class{IA}$ respectively, this interpretation can be lifted from terms to equations. More precisely, let $\btau$ be the $1$-translation of $\mathscr{L}_{\? \class{HA}}$ into $\mathscr{L}_{\? \class{IA}}$ defined as follows for all for $\star \in \{ \land, \lor \}$:
\[
x \star y \longmapsto x \star y \quad \lnot x \longmapsto \Box \lnot x \quad x \to y \longmapsto  \Box  (x \to y)
\]
The interpretation of $\mathcal{IPC}$ into $\mathcal{S}4$ can now be presented as follows:
\begin{equation}\label{Eq:Chap5:Godela}
\Gamma \vdash_{\mathcal{IPC}} \varphi \Longleftrightarrow \sigma\btau_{\ast}(\Gamma) \vdash_{\mathcal{S}4}\sigma\btau_{\ast} (\varphi)
\end{equation}
for every $\Gamma \cup \{ \varphi \} \subseteq Tm( \mathscr{L}_{ \? \class{HA}}, \lambda)$, where $\sigma$ is the substitution sending every variable $x$ to its necessitation $\Box x$. In order to present this translation in our framework, we have to deal with the fact that we allow only translations that send variables to variables. As we mentioned, this  problem is overcome by introducing a \textit{context} in the premises. To explain how, we recall that the terms of $Tm( \mathscr{L}_{ \? \class{HA}}, \lambda)$ are written with variables among $\{ x_{j} : j < \lambda \}$. Then we have that:
\begin{equation}\label{Eq:Chap5:Godelb}
\sigma\btau_{\ast}(\Gamma) \vdash_{\mathcal{S}4}\sigma\btau_{\ast} (\varphi) \Longleftrightarrow \btau_{\ast}(\Gamma) \cup \{ x_{j} \leftrightarrow \Box x_{j} : j < \lambda \}  \vdash_{\mathcal{S}4} \btau_{\ast}(\varphi).
\end{equation}
The left-to-right direction of (\ref{Eq:Chap5:Godelb}) follows from the fact that the algebraic meaning of $x_{j} \leftrightarrow \Box x_{j}$ is $x_{j} \thickapprox \Box x_{j}$. To prove the other direction, suppose that the right-hand deduction holds. Then by structurality we can apply the substitution $\sigma$ to it. This fact, together with $\emptyset \vdash_{\mathcal{S}4} \Box x \leftrightarrow \Box \Box x$, yields the desired conclusion. Now, using the completeness of $\mathcal{IPC}$ and $\mathcal{S}4$ with respect to the corresponding equivalent algebraic semantics, we obtain that
\begin{equation}\label{Eq : 123}
\Phi \vDash_{\class{HA}} \epsilon \thickapprox \delta \Longleftrightarrow \btau^{\ast}(\Phi) \cup  \bigcup_{j < \lambda } \Theta(x_{j})  \vDash_{\class{IA}} \btau^{\ast}(\epsilon \thickapprox \delta ) 
\end{equation}
for every $\Phi \cup \{ \epsilon \thickapprox \delta \} \subseteq Eq (\mathscr{L}_{\?\class{HA}}, \lambda )$, where $\Theta(x) = \{ x \thickapprox \Box x \}$. Observe that (\ref{Eq : 123}) implies condition 1 of Definition \ref{Def : Translation}.  It is easy to see that condition 2, of the same definition, holds as well. Therefore we conclude that $\langle \btau, \Theta \rangle$ is a contextual translation of $\vDash_{\class{HA}}$ into $\vDash_{\class{IA}}$.
\qed
\end{exa}
\begin{exa}[\textsf{Heyting and Boolean Algebras}]\label{Exa : Kol Tr}
The same trick can be applied to subsume Kolmogorov's interpretation of classical propositional calculus $\mathcal{CPC}$ into $\mathcal{IPC}$ \cite{HeKol25} in our framework. Let $\btau$ be the $1$-translation defined as follows:
\[
0 \longmapsto 0 \quad 1 \longmapsto 1 \quad  \lnot x \longmapsto \lnot x \quad x \star y \longmapsto \lnot \lnot (x \star y)
\]
for every $\star \in \{ \land, \lor, \to \}$. The original translation of Kolmogorov states that
\[
\Gamma \vdash_{\mathcal{CPC}} \varphi \Longleftrightarrow \sigma \btau_{\ast}(\Gamma) \vdash_{\mathcal{IPC}}\sigma \btau_{\ast}(\varphi)
\]
for every $\Gamma \cup \{ \varphi \} \subseteq Tm(\mathscr{L}, \lambda)$, where $\sigma$ is the substitution sending every variable $x$ to its double negation $\lnot \lnot x$. Combining it with the observation that $\emptyset \vdash_{\mathcal{IPC}}  \lnot x \leftrightarrow \lnot \lnot \lnot x$, it is easy to see that $\langle \btau, \Theta \rangle$ with $\Theta = \{ x \thickapprox \lnot\lnot x \}$ is a contextual translation of $\vDash_{\class{BA}}$ into $\vDash_{\class{HA}}$, where $\class{BA}$ is the variety of Boolean algebras.
\qed
\end{exa}

The importance of non-trivial contextual $\kappa$-translations of $\vDash_{\class{X}}$ into $\vDash_{\class{Y}}$ is that they correspond to non-trivial right adjoint functors from $\class{Y}$ to $\class{X}$. Notice that right adjoints reverse the direction of contextual translations and vice-versa. We now proceed to establish one half of this correspondence by showing how to construct a right adjoint functor out of a contextual translation. Consider a non-trivial contextual $\kappa$-translation $\langle \btau, \Theta \rangle$ of $\vDash_{\class{X}}$ into $\vDash_{\class{Y}}$. Then consider the set:
\begin{equation}\label{Eq:Sec5.3:language}
\mathscr{L} \coloneqq \{ \btau(\psi) : \psi \in \mathscr{L}_{\?\class{X}} \} \subseteq \mathscr{L}^{\kappa}_{\class{Y}}.
\end{equation}
Observe that $\mathscr{L}$ is a sublanguage of the language of the matrix power $\class{Y}^{[\kappa]}$. Then consider the set 
\[
\theta \coloneqq \{ \vec{\epsilon} \thickapprox \vec{\delta} : \epsilon \thickapprox \delta \in \Theta \}
\]
where $\vec{\epsilon}$ and $\vec{\delta}$ are the $\kappa$-sequences constantly equal to $\epsilon$ and $\delta$ respectively. Observe that $\theta$ is a set of identities between $\kappa$-sequences of terms of $\class{Y}$ in $\kappa$ variables. Now, $\kappa$-sequences of terms of $\class{Y}$ in $\kappa$-many variables can be viewed as \textit{unary} terms of the matrix power $\class{Y}^{[\kappa]}$. Thus $\theta$ can be viewed as a set of equations in one variable in the language of $\class{Y}^{[\kappa]}$. Hence we have the three basic ingredients of our construction: a matrix power $\class{Y}^{[\kappa]}$, a sublanguage $\mathscr{L} \subseteq \mathscr{L}^{\kappa}_{\class{Y}}$, and a set of equations $\theta \subseteq Eq(\mathscr{L}^{\kappa}_{\class{Y}}, 1)$.

There is still a technicality we must take into account: when $\kappa$ is infinite the matrix power $\class{Y}^{[\kappa]}$ may fail to be a generalized quasi-variety. Let $\class{K}$ be the class of algebras defined as follows:
\[
\class{K} \coloneqq \left\{ \begin{array}{ll}
 \QQQ(\class{Y}^{[\kappa]}) & \text{if $\class{X}$ and $\class{Y}$ are quasi-varieties}\\
 & \text{and $\textup{Cg}^{\Tm_{\class{Y}}(\kappa)}_{\class{Y}}(\Theta)$ is finitely generated}\\[1ex]
 \GGG\QQQ_{\lambda}(\class{Y}^{[\kappa]}) & \text{otherwise, where }\lambda\text{ is infinite and }\UUU_{\lambda}(\class{X}) = \class{X}
  \end{array} \right.  
\]
where the expressions $\QQQ$ and $\GGG\QQQ_{\lambda}$ have been introduced at pag.\ \pageref{DefOfQQQ}. Observe that in the above definition $\lambda$ is not uniquely determined, but any choice will be equivalent for our purposes.

\begin{Theorem}\label{Thm : Tr > RA}
Let $\class{X}$ and $\class{Y}$ be generalized quasi-varieties, let $\langle \btau, \Theta \rangle$ be a non-trivial contextual $\kappa$-translation of $\vDash_{\class{X}}$ into $\vDash_{\class{Y}}$, and let $\class{K}$ be the class just introduced. The maps $[\kappa] \colon \class{Y} \to \class{K}$ and $\theta_{\mathscr{L}} \colon \class{K} \to \class{X}$ defined above are right adjoint functors. In particular, the composition $\theta_{\mathscr{L}} \circ [\kappa] \colon \class{Y} \to \class{X}$ is a non-trivial right adjoint.
\end{Theorem}
\begin{proof}
Observe that $\class{K}$ is a generalized quasi-variety. Therefore we can apply Theorem \ref{Thm:McKenzie}, yielding that $[\kappa] \colon \class{Y} \to \class{K}$ is a right adjoint functor. Now we turn to prove the same for $\theta_{\mathscr{L}}$. We will detail the case where $\class{X}$ and $\class{Y}$ are quasi-varieties and $\textup{Cg}^{\Tm_{\class{Y}}(\kappa)}_{\class{Y}}(\Theta)$ finitely generated, since the other case is analogous. Since $\class{Y}$ is a quasi-variety and $\textup{Cg}^{\Tm_{\class{Y}}(\kappa)}_{\class{Y}}(\Theta)$ is finitely generated, there is a finite set $\{ \langle \alpha_{i}, \beta_{i} \rangle : i < n \} \subseteq \Theta$ such that $\{ \langle \alpha_{i}, \beta_{i} \rangle : i < n \} \semeq_{\class{Y}} \Theta$. It is easy to see that
\begin{equation}\label{Eq : Fin theta}
\{  \vec{\alpha}_{i}  \thickapprox  \vec{\beta}_{i}  : i < n \} \semeq_{\class{Y}^{[\kappa]}} \theta
\end{equation}
where $\vec{\alpha}_{i}$ and $\vec{\beta}_{i}$ are the $\kappa$-sequences constantly equal to $\alpha_{i}$ and $\beta_{i}$ respectively.

Now from condition 2 of Definition \ref{Def : Translation} it follows that $\theta$ is compatible with $\mathscr{L}$ in $\class{Y}^{[\kappa]}$, where $\mathscr{L}$ is the language defined in (\ref{Eq:Sec5.3:language}). From (\ref{Eq : Fin theta}) we know that this compatibility condition can be expressed by a set of deductions, whose antecedent is finite, of the equational consequence relative to $\class{Y}^{[\kappa]}$, i.e.,
\[
\bigcup_{ j \leq m}\{ \vec{\alpha}_{i}  \thickapprox  \vec{\beta}_{i}  : i < n \}(\vec{x}_{j}) \vDash_{\class{Y}^{[\kappa]}} \theta(\btau(\psi)(\vec{x}_{1}, \dots, \vec{x}_{n}))
\]
for every $m$-ary $\psi \in \mathscr{L}$. In particular, this implies that $\theta$ is still compatible with $\mathscr{L}$ in $\class{K}$ (recall that $\class{K}$ is the \textit{quasi-variety} generated by $\class{Y}^{[\kappa]}$).

We claim that $\theta_{\mathscr{L}}(\A) \in \class{X}$ for every $\A \in \class{K}$. To prove this, consider any finite deduction
\[
\varphi_{1} \thickapprox \psi_{1}, \dots, \varphi_{m} \thickapprox \psi_{m} \vDash_{\class{X}} \epsilon \thickapprox \delta.
\]
Let $x_{1}, \dots, x_{p}$ be the variables that occur in it. From condition 1 of Definition \ref{Def : Translation} it follows that
\[
\{ \btau_{\ast}(\varphi_{t}) \thickapprox \btau_{\ast}(\psi_{t}) : t \leq m \} \cup \bigcup_{ j \leq p} \theta (\vec{x}_{j})   \vDash_{\class{Y}^{[\kappa]}} \btau_{\ast}(\epsilon) \thickapprox \btau_{\ast}(\delta)
\]
where $\vec{x}_{j} = \langle x^{i}_{j} : i < \kappa \rangle$. Thanks to (\ref{Eq : Fin theta}) the above deduction can be expressed by a collection of  deductions, whose antecedent is finite, of the equational consequence relative to $\class{Y}^{[\kappa]}$, i.e.,
\[
\{ \btau_{\ast}(\varphi_{t}) \thickapprox \btau_{\ast}(\psi_{t}) : t \leq m \} \cup \bigcup_{ j \leq p} \{ \vec{\alpha}_{i}  \thickapprox  \vec{\beta}_{i}  : i < n \}(\vec{x}_{j})  \vDash_{\class{Y}^{[\kappa]}} \btau_{\ast}(\epsilon) \thickapprox \btau_{\ast}(\delta).
\]
Since $\class{K}$ is the \textit{quasi-variety} generated by $\class{Y}^{[\kappa]}$, we know that the above deduction persists in $\class{K}$. Together with the fact that $\{ \vec{\alpha}_{i}  \thickapprox  \vec{\beta}_{i}  : i < n \} \subseteq \theta$, this implies that for every $\A \in \class{K}$ and every $a_{1}, \dots, a_{p} \in \theta_{\mathscr{L}}(A)$, we have that:
\begin{align*}
\text{if }&\theta_{\mathscr{L}}(\A)\vDash \varphi_{1} \thickapprox \psi_{1}, \dots, \varphi_{m} \thickapprox \psi_{m} \llbracket a_{1}, \dots, a_{p} \rrbracket,\\
\text{then }& \theta_{\mathscr{L}}(\A) \vDash \epsilon \thickapprox \delta\llbracket a_{1}, \dots, a_{p} \rrbracket.
\end{align*}
Thus we showed that $\theta_{\mathscr{L}}(\A)$ satisfies every quasi-equation that holds in $\class{X}$. Since $\class{X}$ is a quasi-variety, we conclude that $\theta_{\mathscr{L}}(\A) \in \class{X}$. This establishes our claim. Hence we can apply Theorem \ref{Thm:ThetaRight}, yielding that
$\theta_{\mathscr{L}} \colon \class{K} \to \class{X}$ is a right adjoint functor. We conclude that $\theta_{\mathscr{L}} \circ [\kappa] \colon \class{Y} \to \class{X}$ is a right adjoint functor.

The fact that $\theta_{\mathscr{L}} \circ [\kappa]$ is non-trivial follows from the fact that so is $\langle \btau, \Theta \rangle$.
\end{proof}

If we apply the above construction to G\"odel and Kolmogorov's translations, we obtain some well-known transformations:

\begin{exa}[\textsf{Open and Regular Elements}]\label{Exa : Open Elements}
Given $\A \in \class{IA}$, an element $a \in A$ is \textit{open} if $\Box a = a$. The set of open elements $\textup{Op}(\A)$ of $\A$ is closed under the lattice operations and contains the bounds. Moreover we can equip it with an implication $\multimap$ and with a negation $\thicksim$ defined for every $a, b \in \textup{Op}(\A)$ as follows:
\[
a \multimap b \coloneqq \Box^{\A}( a \to^{\A} b) \text{ and }\thicksim a \coloneqq \Box^{\A} \lnot^{\A} a.
\]
It is well known that
\[
\textup{Op}(\A) \coloneqq \langle \textup{Op}(\A), \land, \lor, \multimap, \thicksim, 0, 1 \rangle
\]
is a Heyting algebra. Now, every homomorphism $f \colon \A \to \B$ between interior algebras restricts to a homomorphism $f \colon \textup{Op}(\A) \to \textup{Op}(\B)$. Therefore the map $\textup{Op} \colon \class{IA} \to \class{HA}$ can be regarded as a functor. As the reader may have guessed, it is in fact the right adjoint functor induced by G\"odel's translation of $\mathcal{IPC}$ into $\mathcal{S}4$ (Example \ref{Exa : God Tr}).

A similar correspondence arises from Kolmogorov's translation of $\mathcal{CPC}$ into $\mathcal{IPC}$. More precisely, given $\A \in \class{HA}$, an element $a \in A$ is \textit{regular} if $\lnot \lnot a = a$. It is well known that the set of regular elements $\textup{Reg}(\A)$ of $\A$ is closed under $\land, \lnot$ and $\to$ and contains the bounds. Moreover we can equip it with a new join $\sqcup$ defined for every $a, b \in \textup{Reg}(\A)$ as follows:
\[
a \sqcup b \coloneqq \lnot^{\A} \lnot^{\A}( a \lor b).
\]
It is well known that
\[
\textup{Reg}(\A) \coloneqq \langle \textup{Reg}(\A), \land, \sqcup, \to, \lnot, 0, 1 \rangle
\]
is a Boolean algebra. Now, every homomorphism $f \colon \A \to \B$ between Heyting algebras restricts to a homomorphism $f \colon \textup{Reg}(\A) \to \textup{Reg}(\B)$. Therefore the map $\textup{Reg} \colon \class{HA} \to \class{BA}$ can be regarded as a functor, which is exactly the right adjoint functor induced by Kolmogorov's translation (Example \ref{Exa : Kol Tr}).
\qed
\end{exa}

\section{From right adjoints to translations}\label{Sec : RA}

In this section we will work with a fixed non-trivial left adjoint functor $\mathcal{F} \colon \class{X} \to \class{Y}$ between generalized quasi-varieties. Our goal is to construct a contextual translation of $\vDash_{\class{X}}$ into $\vDash_{\class{Y}}$ induced by $\mathcal{F}$. We rely on the following easy observation:

\begin{Lemma}\label{Lem : Not empty}
Let $\mathcal{F} \colon \class{X} \to \class{Y}$ be a non-trivial left adjoint functor between generalized quasi-varieties. The universe of $\mathcal{F}(\Tm_{\class{X}}(1))$ is non-empty.
\end{Lemma}

Now we construct the contextual translation $\langle \btau, \Theta \rangle$ induced by $\mathcal{F} \colon \class{X} \to \class{Y}$. By Lemma \ref{Lem : Not empty} we know that $\mathcal{F}(\Tm_{\class{X}}(1)) \ne \emptyset$. Then we can choose a cardinal $\kappa > 0$ and a surjective homomorphism $\pi_{1} \colon  \Tm_{\class{Y}}(\kappa) \to \mathcal{F}(\Tm_{\class{X}}(1))$. Let $\Theta$ be the kernel of $\pi$ and observe that it can be viewed as a set of equations in $Eq(\mathscr{L}_{ \? \class{Y}}, \kappa)$.

In order to construct the $\kappa$-translation $\btau$ of $\mathscr{L}_{\?\class{X}}$ into $\mathscr{L}_{\?\class{Y}}$, consider a cardinal $\lambda > 0$. Since $\mathcal{F}$ preserves copowers and the algebra $\Tm_{\class{X}}(\lambda)$ is the $\lambda$-th copower of $\Tm_{\class{X}}(1)$, we know that $\mathcal{F}(\Tm_{\class{X}}(\lambda))$ is the $\lambda$-th copower of $\mathcal{F}(\Tm_{\class{X}}(1))$. Keeping in mind how coproducts look like in prevarieties, we can identify $\mathcal{F}(\Tm_{\class{X}}(\lambda))$ with the quotient of the free algebra $\Tm_{\class{Y}}(\kappa \times \lambda)$ with free generators $\{ x^{i}_{j} : i < \kappa, j < \lambda \}$ under the $\class{Y}$-congruence generated by $\bigcup_{j  < \lambda} \Theta(\vec{x}_{j})$ where $\vec{x}_{j} = \langle x_{j}^{i} : i < \kappa \rangle$.

The above construction can be carried out also for $\lambda = 0$ as follows. Recall that $\mathcal{F}$ preserves initial objects, since these are special colimits. Thus we can assume that $\mathcal{F}(\Tm_{\class{X}}(0)) = \Tm_{\class{Y}}(0)$. Now we have that $\Tm_{\class{Y}}(0)$ is exactly the quotient of $\Tm_{\class{Y}}(\kappa \times 0 )$ under the $\class{Y}$-congruence generated by the union of zero-many copies of $\Theta$, i.e., under the identity relation. Thus we identify $\mathcal{F}(\Tm_{\class{X}}(\lambda))$ with a quotient of $\Tm_{\class{Y}}( \kappa\times\lambda  )$ for every cardinal $\lambda$. Accordingly, we denote by $\pi_{\lambda} \colon \Tm_{\class{Y}}(\kappa\times\lambda) \to \mathcal{F}(\Tm_{\class{X}}(\lambda))$ the corresponding canonical map.

\begin{Definition}\label{Def:Chap5:ArrowsFunctions}
Let $\lambda$ be a cardinal and $\varphi \in Tm(\mathscr{L}_{\class{X}}, \lambda)$. We denote also by $\varphi \colon \Tm_{\class{X}}(1) \to \Tm_{\class{X}}(\lambda)$ the unique homomorphism that sends $x$ to $\varphi$, where $x$ is the free generator of $\Tm_{\class{X}}(1)$.
\end{Definition}

Now, we are ready to construct the $\kappa$-translation $\btau$ of $\mathscr{L}_{\?\class{X}}$ into $\mathscr{L}_{\?\class{Y}}$. Consider an $n$-ary basic operation $\psi \in \mathscr{L}_{ \? \class{X}}$. Since $\pi_{n}$ is surjective and $\Tm_{\class{Y}}(\kappa)$ is onto-projective in $\class{Y}$, there is a homomorphism
\[
\btau(\psi)\colon \Tm_{\class{Y}}(\kappa) \to \Tm_{\class{Y}}(\kappa \times n)
\]
that makes the following diagram commute:
\begin{equation}\label{Eq:Chap5:xypic}
\xymatrix@R=35pt @C=45pt @!0{
&&\Tm_{\class{Y}}(\kappa)\ar[d]^{\pi_{1}}\ar@{.>}@/^-2pc/[ddll]|-{\btau(\psi)}\\
&& \mathcal{F}(\Tm_{\class{X}}(1)) \ar[d]^{\mathcal{F}(\psi)}\\
\Tm_{\class{Y}}(\kappa \times n)\ar[rr]_{\pi_{n}}&& \mathcal{F}(\Tm_{\class{X}}(n))
}
\end{equation}

\noindent The map $\btau(\psi) $ can be identified with its values on the generators $\{ x^{i} : i < \kappa \}$ of $\Tm_{\class{Y}}(\kappa)$. In this way it becomes a $\kappa$-sequence
\[
\langle \btau(\psi) (x^{i}) : i < \kappa \rangle
\]
of terms in variables $\{ x_{j}^{i} : i< \kappa, 1 \leq j \leq n \}$.

Let $\btau$ be the $\kappa$-translation of $\mathscr{L}_{ \? \class{X}}$ into $\mathscr{L}_{ \? \class{Y}}$ obtained by applying this construction to every $\psi \in \mathscr{L}_{ \? \class{X}}$. Hence we constructed a pair $\langle \btau, \Theta\rangle$, where $\btau$ is a $\kappa$-translation of $\mathscr{L}_{ \? \class{X}}$ into $\mathscr{L}_{ \? \class{Y}}$ and $\Theta \subseteq Eq(\mathscr{L}_{ \? \class{Y}}, \kappa)$.

\begin{Theorem}\label{Lemma : RA > TR}
Let $\mathcal{F} \colon \class{X} \to \class{Y}$ be a non-trivial left adjoint functor between generalized quasi-varieties. The pair $\langle \btau, \Theta \rangle$ defined above is a non-trivial contextual translation of $\vDash_{\class{X}}$ into $\vDash_{\class{Y}}$.
\end{Theorem}

\begin{proof}[Proof sketch.]
Consider a cardinal $\lambda$. We know that $\btau$ can be extended to a function $\btau_{\ast} \colon Tm(\mathscr{L}_{\?\class{X}}, \lambda) \to Tm(\mathscr{L}_{\?\class{Y}}, \kappa \times \lambda)^{\kappa}$, where the terms $Tm(\mathscr{L}_{\?\class{X}}, \lambda)$ and $Tm(\mathscr{L}_{\?\class{Y}}, \kappa \times \lambda)$ are built respectively with variables among $\{ x_{j} : j < \lambda \}$ and $\{ x_{j}^{i} : i< \kappa, j < \lambda \}$.

Now, consider $\varphi \in Tm(\mathscr{L}_{\?\class{X}}, \lambda)$. Observe that $\btau_{\ast}(\varphi)$ is a $\kappa$-sequence of terms of $\class{Y}$ in variables $\{ x_{j}^{i} : i< \kappa, j < \lambda \}$. Thus $\btau_{\ast}(\varphi)$ can be regarded as a map from the free generators of $\Tm_{\class{Y}}(\kappa)$ to $\Tm_{\class{Y}}( \kappa\times\lambda)$. Since $\Tm_{\class{Y}}(\kappa)$ is a free algebra, this assignment extends uniquely to a homomorphism
\[
\btau_{\ast}(\varphi) \colon \Tm_{\class{Y}}(\kappa) \to \Tm_{\class{Y}}(\kappa\times \lambda).
\]
Recall that $\mathcal{F}$ preserves colimits, since it is left adjoint. Keeping this in mind, it is not difficult to prove the following:

\begin{fact}\label{Claim:1c}
For every cardinal $\lambda$ and every $\varphi \in Tm(\mathscr{L}_{\?\class{X}}, \lambda)$, the following diagram commutes:
\[
\xymatrix@R=45pt @C=95pt @!0{
\Tm_{\class{Y}}(\kappa) \ar[r]^{\btau_{\ast}(\varphi)} \ar[d]_{\pi_{1}} & \Tm_{\class{Y}}(\kappa\times\lambda)\ar[d]^{\pi_{\lambda}}\\
\mathcal{F}(\Tm_{\class{X}}(1)) \ar[r]_{\mathcal{F}(\varphi)} & \mathcal{F}(\Tm_{\class{X}}(\lambda))
}
\]
\end{fact}

Now we turn to prove that $\langle \btau, \Theta\rangle$ is a contextual translation of $\vDash_{\class{X}}$ into $\vDash_{\class{Y}}$. We begin by  showing that $\langle \btau, \Theta\rangle$ satisfies 1 of Definition \ref{Def : Translation}. To this end, consider a cardinal $\lambda$ and equations $\Phi \cup \{ \epsilon \thickapprox \delta \} \subseteq Eq( \mathscr{L}_{\?\class{X}}, \lambda)$ such that $\Phi \vDash_{\class{X}} \epsilon \thickapprox \delta$. Define $\mu \coloneqq \vert \Phi \vert$. For the sake of simplicity we identify $\mu$ with the set $\Phi$. Then consider the map $\btau_{\ast} \colon Tm(\mathscr{L}_{\?\class{X}}, \lambda) \to Tm(\mathscr{L}_{\?\class{Y}}, \kappa\times\lambda)^{\kappa}$. Consider also the free algebras $\Tm_{\class{X}}(\mu)$ and $\Tm_{\class{Y}}(\kappa \times \mu)$ with free generators $\{ x_{\alpha\thickapprox \beta} : \alpha \thickapprox \beta \in \Phi \}$ and $\{ x_{\alpha\thickapprox \beta}^{i} : i < \kappa, \alpha \thickapprox \beta \in \Phi \}$ respectively. Then let
\[
p_{l}, p_{r} \colon \Tm_{\class{X}}(\mu) \rightrightarrows \Tm_{\class{X}}(\lambda) \?\?\?\?\?\? \text{ and }\?\?\?\?\?\? q_{l}, q_{r} \colon \Tm_{\class{Y}}(\kappa \times \mu) \rightrightarrows \Tm_{\class{Y}}(\kappa\times\lambda )
\]
be the homomorphisms defined respectively by the following rules:
\begin{align*}
p_{l}(x_{\alpha \thickapprox \beta}) \coloneqq \alpha &\text{ and }q_{l}(x_{\alpha \thickapprox \beta}^{i}) \coloneqq \btau_{\ast}(\alpha)(i)\\
p_{r}(x_{\alpha \thickapprox \beta}) \coloneqq \beta &\text{ and }q_{r}(x_{\alpha \thickapprox \beta}^{i}) \coloneqq \btau_{\ast}(\beta)(i).
\end{align*}
Observe that 
\begin{equation}\label{Eq:Chap5:MapsPQ}
\pi_{\lambda} \circ q_{l} = \mathcal{F}(p_{l}) \circ \pi_{\mu} \?\?\?\? \text{ and } \?\?\?\? \pi_{\lambda} \circ q_{r} = \mathcal{F}(p_{r}) \circ \pi_{\mu}.
\end{equation}

Now, let $\phi$ be the $\class{X}$-congruence of $\Tm_{\class{X}}(\lambda)$ generated by $\Phi$. It is clear that $\pi_{\phi}$ is a coequalizer of $p_{l}$ and $p_{r}$. Since $\mathcal{F}$ preserves colimits and $\pi_{\mu}$ is surjective, this means that $\mathcal{F}(\pi_{\phi})$ is also a coequalizer of $\mathcal{F}(p_{l}) \circ \pi_{\mu}$ and $\mathcal{F}(p_{r}) \circ \pi_{\mu}$. Finally, with an application of (\ref{Eq:Chap5:MapsPQ}), we conclude that $\mathcal{F}(\pi_{\phi})$ is a coequalizer of $\pi_{\lambda} \circ q_{l}$ and $\pi_{\lambda} \circ q_{r}$. In particular, this implies that the kernel of $\mathcal{F}(\pi_{\phi}) \circ \pi_{\lambda}$ is the $\class{Y}$-congruence of $\Tm_{\class{Y}}(\kappa\times\lambda)$ generated by
\begin{equation}\label{Eq:Chap5:genconproof}
\btau^{\ast}(\Phi) \cup \bigcup_{j  < \lambda} \Theta(\vec{x}_{j})
\end{equation}
where $\vec{x}_{j} = \langle x_{j}^{i} : i < \kappa \rangle$. Now, recall that $\Phi \vDash_{\class{X}} \epsilon \thickapprox \delta$ and, therefore, that $\langle \epsilon, \delta \rangle \in \phi$. This means that $\pi_{\phi} \circ \epsilon = \pi_{\phi} \circ \delta$, where $\epsilon, \delta \colon \Tm_{\class{X}}(1) \rightrightarrows \Tm_{\class{X}}(\lambda)$. By Fact \ref{Claim:1c} this implies that
\[
\mathcal{F}(\pi_{\phi}) \circ \pi_{\lambda} \circ \btau_{\ast}(\epsilon) =  \mathcal{F}(\pi_{\phi}) \circ \pi_{\lambda} \circ \btau_{\ast}(\delta). 
\]
Together with the description of the kernel of $\mathcal{F}(\pi_{\phi}) \circ \pi_{\lambda}$ given in (\ref{Eq:Chap5:genconproof}), this yields
\[
\btau^{\ast}(\Phi) \cup \bigcup_{j  < \lambda} \Theta(\vec{x}_{j}) \vDash_{\class{Y}} \btau^{\ast}( \epsilon \thickapprox \delta).
\]
Hence $\langle \btau, \Theta\rangle$ satisfies 1 of Definition \ref{Def : Translation}. 

%The proof of the fact that it satisfies condition 2, of the same definition, is an easy exercise.

To prove that $\langle \btau, \Theta\rangle$ satisfies condition 2 of the same definition, consider an $n$-ary operation symbol $\psi \in \mathscr{L}_{\? \class{X}}$ and $\epsilon \thickapprox \delta \in \Theta$. Fact \ref{Claim:1c} and the observation that the kernel of $\pi_{1}$ is the $\class{Y}$-congruence of $\Tm_{\class{Y}}(\kappa)$ generated by $\Theta$ imply that
\[
\pi_{n} ( \epsilon (\btau_{\ast}(\psi) / \vec{x} )) = \pi_{n} ( \delta ( \btau_{\ast}(\psi) / \vec{x})).
\]
Since $\pi_{n}$ is the kernel of the $\class{Y}$-congruence of $\Tm_{\class{Y}}(\kappa \times n)$ generated  by $\Theta(\vec{x}_{1}) \cup \cdots \cup \Theta(\vec{x}_{n})$, we conclude that
\[
\Theta(\vec{x}_{1}) \cup \cdots \cup \Theta(\vec{x}_{n}) \vDash_{\class{Y}} \epsilon (\btau_{\ast}(\psi) / \vec{x}) \thickapprox \delta (\btau_{\ast}(\psi) / \vec{x} ).
\]
This establishes that $\langle \btau, \Theta\rangle$ is a contextual translation of $\vDash_{\class{X}}$ into $\vDash_{\class{Y}}$.

It only remains to prove that $\langle \btau, \Theta\rangle$ is non-trivial. But this is a consequence of the fact that $\mathcal{F}$ is non-trivial.
\end{proof}

As an exemplification of the construction above, we will describe the contextual translation associated with the adjunction between Kleene algebras and bounded distributive lattices.

\begin{exa}[\textsf{Kleene Algebras}]\label{Exa : KA2}
Let $\mathcal{G} \colon \class{DL}_{01} \to \class{KA}$ be the functor described in Example \ref{Exa : KA}. In \cite{Ci86} a functor $\mathcal{F}$ left adjoint to $\mathcal{G}$ is described. Let us briefly recall its behaviour. Given $\A \in \class{KA}$, we let $\text{Pr}(\A)$ be the Priestley space dual to the bounded lattice reduct of $\A$ \cite{DaPr02}. Moreover, we equip it with a map $g \colon \text{Pr}(\A) \to \text{Pr}(\A)$ defined by the rule
\[
g(F) \longmapsto A \smallsetminus \{ \lnot a : a \in F \} \text{, with } F \in  \text{Pr}(\A).
\]
Now observe that
\[
\text{Pr}(\A)^{+} \coloneqq \{ F \in \text{Pr}(\A) : F \subseteq g(F) \}
\]
is the universe of a Priestley subspace of $\text{Pr}(\A)$. Keeping this in mind, we let $\mathcal{F}(\A)$ be the bounded distributive lattice dual to $\text{Pr}(\A)^{+}$. Moreover, given a homomorphism $f \colon \A \to \B$ in $\class{KA}$, we let $\mathcal{F}(f) \colon \mathcal{F}(\A) \to  \mathcal{F}(\B)$ be the map defined by the rule
\[
U \longmapsto \{ F \in \text{Pr}(\B)^{+} : f^{-1}(F) \in U \} \text{, for each }U \in \mathcal{F}(\A).
\]
The map $\mathcal{F} \colon \class{KA} \to \class{DL}_{01}$ is the functor left adjoint to $\mathcal{G}$.

Now we turn to describe the contextual translation associated with the adjunction $\mathcal{F} \dashv \mathcal{G}$. To this end, observe that the free Kleene algebra $\Tm_{\class{KA}}(1)$, its image $\mathcal{F}(\Tm_{\class{KA}}(1))$ in $\class{DL}_{01}$ and the free bounded distributive lattice $\Tm_{\class{DL}_{01}}(2)$ are respectively the algebras depicted below.

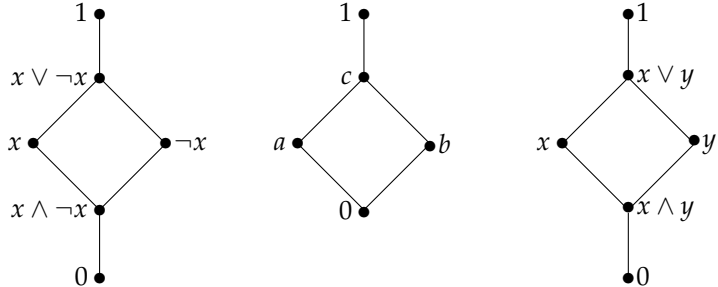
\[
\xymatrix@R=25pt @C=25pt @!0{
&  *-{\mathllap{\ \ 1\,}\bullet  }\ar@{-}[d]& &&&*-{\mathllap{\ \ 1\,}\bullet  }\ar@{-}[d]&&&&*-{\mathrlap{\ \  1\,}\bullet  }\ar@{-}[d]&&\\
&  *-{\mathllap{\ \ x \lor \lnot x\,}\bullet  }\ar@{-}[dr]\ar@{-}[dl]&&&&*-{\mathllap{\ \ c\,}\bullet  }\ar@{-}[dr]\ar@{-}[dl]&&&&*-{\mathrlap{\ \   x\lor y \,}\bullet  }\ar@{-}[dl]\ar@{-}[dr]&&\\
*-{\mathllap{\ \ x\,}\bullet  } \ar@{-}[dr]&  &*-{\mathrlap{\ \ \lnot x\,}\bullet  }\ar@{-}[dl]&&*-{\mathllap{\ \ a\,}\bullet  }\ar@{-}[dr]&&*-{\mathrlap{\ \ b\,}\bullet  }\ar@{-}[dl]&&*-{\mathllap{\ \  x\,}\bullet  }\ar@{-}[dr]&&*-{\mathrlap{\ \  y\,}\bullet  }\ar@{-}[dl]&\\
& *-{\mathllap{\ \ x \land \lnot x\,}\bullet  } \ar@{-}[d]&&&&*-{\mathllap{\ \ 0\,}\bullet  }&&&&*-{\mathrlap{\ \  x\land y\,}\bullet  }\ar@{-}[d]&&\\
& *-{\mathllap{\ \ 0\,}\bullet  }&&&&&&&&*-{\mathrlap{\ \  0\,}\bullet  }&&
}
\]
Then let $\pi \colon \Tm_{\class{DL}_{01}}(2) \to \mathcal{F}(\Tm_{\class{KA}}(1))$ be the unique (surjective) homomorphism determined by the assignment $\pi(x) = a$ and $\pi(y)= b$. Following the general construction described above, we should identify $\Theta$ with the kernel of $\pi$ viewed as a set of 
equations in 2 variables. But the only equation of this kind that is not vacuously satisfied is $x \land y \thickapprox 0$. Hence we can set without loss of generality $\Theta \coloneqq \{ x \land y \thickapprox 0 \}$.

The description of $\btau$ is more complicated and we will detail it only for the case of negation. First 
observe that $\lnot \colon \Tm_{\class{KA}}(1) \to \Tm_{\class{KA}}(1)$ is the unique endomorphism that sends $x$ to $\lnot x$. Then, applying the 
definition of $\mathcal{F}$, it is easy to see that $\mathcal{F}(\lnot)$ is the endomorphism of $\mathcal{F}(\Tm_{\class{KA}}(1))$ that behaves as the identity except that it interchanges $a$ and $b$. Now we have to choose an endomorphism $\btau(\lnot)$ of  $\Tm_{\class{DL}_{01}}(2)$ such that $\pi \circ \btau(\lnot) = \mathcal{F}(\lnot)\circ \pi$. It is easy to see that the unique homomorphism $\btau(\lnot)$ determined by the assignment $\btau(\lnot)(x) = y$ and $\btau(\lnot)(y) = x$ fulfils this condition. Hence the translation of $\lnot$ consists of the pair $\langle y, x\rangle$. The same idea allows us to extend $\btau$ to the other constant and binary basic symbols of $\class{KA}$ as follows:\footnote{At this stage the reader may find it useful to compare the translation displayed here with the sublanguage $\mathscr{L}$ of the matrix power $\class{DL}_{01}$ that we considered in Example \ref{Exa : KA}.}
\begin{alignat*}{3}
x \sqcap y &\longmapsto \langle x^{1}, x^{2}\rangle \sqcap \langle y^{1}, y^{2}\rangle &\coloneqq  &\langle x^{1} \land y^{1}, x^{2} \lor y^{2}\rangle\\
x \sqcup y &\longmapsto  \langle x^{1}, x^{2}\rangle \sqcup \langle y^{1}, y^{2}\rangle &\coloneqq &\langle x^{1} \lor y^{1}, x^{2} \land y^{2}\rangle
\end{alignat*}
and
\[
\lnot x \mapsto \lnot \langle x^{1}, x^{2}\rangle \coloneqq  \langle x^{2}, x^{1}\rangle \qquad 1 \mapsto 1 \coloneqq \langle 1, 0\rangle \qquad 0 \mapsto 0 \coloneqq \langle 0, 1 \rangle.
\]
By Theorem \ref{Lemma : RA > TR} the pair $\langle \btau, \Theta \rangle$ is a contextual translation of $\vDash_{\class{KA}}$ into $\vDash_{\class{DL}_{01}}$. 
%
%For the reader familiar with the theory of algebraizable logics \cite{BP89} it may be interesting to observe that this contextual translation is not induced by a translation between two propositional logics (as was the case in Examples \ref{Exa : God Tr} and \ref{Exa : Kol Tr}). This is due to the fact that $\class{DL}_{01}$ and $\class{KA}$ are not the equivalent algebraic semantics of any algebraizable logics.
\qed
\end{exa}

\section{Decomposition of right adjoints}\label{Sec : Main}

In the preceding sections we drew a correspondence between adjunctions and contextual translations. Now we are ready to present the main outcome of this correspondence, namely the observation that every every right adjoint functor between generalized quasi-varieties can be decomposed  into a combination of two canonical deformations, i.e\ the matrix power with (possibly) infinite exponents and the $\theta_{\mathscr{L}}$ construction:
\begin{Theorem}\label{Thm:Main}
Let $\class{X}$ and $\class{Y}$ be generalized quasi-varieties.
\benormal

\item For every non-trivial right adjoint $\mathcal{G} \colon \class{Y} \to \class{X}$ there are a generalized quasi-variety $\class{K}$ and functors $[\kappa] \colon \class{Y} \to \class{K}$ and $\theta_{\mathscr{L}} \colon \class{K} \to \class{X}$ (where $\theta$ is compatible with $\mathscr{L}$ in $\class{K})$ such that $\mathcal{G}$ is naturally isomorphic to $\theta_{\mathscr{L}}\circ [\kappa]$.

\item Every functor of the form $\theta_{\mathscr{L}} \circ [\kappa] \colon \class{Y} \to \class{X}$ (where $\theta$ is compatible with $\mathscr{L}$ in $\class{Y}^{[\kappa]}$) is a right adjoint.
\enormal
\end{Theorem}
\begin{proof}
1. Let $\mathcal{F}$ be the functor left adjoint to $\mathcal{G}$ and let $\eta$, $\epsilon$ be the unit and counit of the adjunction respectively. In Theorem \ref{Lemma : RA > TR} we showed that $\mathcal{F}$ gives rise to a contextual translation $\langle \btau, \Theta \rangle$ of $\vDash_{\class{X}}$ into $\vDash_{\class{Y}}$.  Then consider the generalized quasi-variety $\class{K}$ and the right adjoint functors $[\kappa] \colon \class{Y} \to \class{K}$ and $\theta_{\mathscr{L}} \colon \class{K} \to \class{X}$ associated with $\langle \btau, \Theta \rangle$ as in Theorem \ref{Thm : Tr > RA}. We will prove that $\mathcal{G}$ and the composition $\theta_{\mathscr{L}} \circ [\kappa]$ are naturally isomorphic.

To this end, it will be convenient to work with some substitutes of $\mathcal{G}$ and $\theta_{\mathscr{L}} \circ [\kappa]$. Let $\class{ALG}_{\class{X}}$ be the category of all algebras of the type of $\class{X}$. Then let $\mathcal{G}^{\ast} \colon \class{Y} \to \class{ALG}_{\class{X}}$ be the functor defined by the rule
\begin{align*}
\A & \longmapsto \hom(\Tm_{\class{X}}(1), \mathcal{G}(\A))\\
f & \longmapsto \mathcal{G}(f) \circ (\cdot)
\end{align*}
for every algebra $\A$ and homomorphism $f$ in $\class{Y}$. The operations of the algebra $\mathcal{G}^{\ast}(\A)$ are defined as follows. Given an $n$-ary operation $\psi \in \mathscr{L}_{\? \class{X}}$ with corresponding arrow $\psi \colon \Tm_{\class{X}}(1) \to \Tm_{\class{X}}(n)$, we set
\[
\psi^{\mathcal{G}^{\ast}(\A)}(f_{1}, \dots, f_{n}) \coloneqq \langle f_{1}, \dots, f_{n}\rangle \circ \psi 
\]
for every $f_{1}, \dots, f_{n} \in \mathcal{G}^{\ast}(\A)$, where $\langle f_{1}, \dots, f_{n}\rangle \colon \Tm_{\class{X}}(n) \to \mathcal{G}(\A)$ is the map induced by the universal property coproduct. 

Now observe that the map $\zeta_{\A} \colon \mathcal{G}(\A) \to \mathcal{G}^{\ast}(\A)$ that takes an element $a \in \mathcal{G}(\A)$ to the unique arrow $f \in \mathcal{G}^{\ast}(\A)$ such that $f(x) = a$ is an isomorphism for every $\A \in \class{Y}$. It is easy to see that the global map $\zeta \colon \mathcal{G} \to \mathcal{G}^{\ast}$ is a natural isomorphism between $\mathcal{G}, \mathcal{G}^{\ast} \colon \class{Y} \to \class{ALG}_{\class{X}}$. As a consequence, we obtain the following:
\begin{fact}\label{Fact:Chap5:natiso1}
The map $\mathcal{G}^{\ast}$ can be viewed as a functor from $\class{Y}$ to $\class{X}$ naturally isomorphic to $\mathcal{G}$.
\end{fact}

Then we construct our substitute for $\theta_{\mathscr{L}} \circ [\kappa]$. Consider the functor
\[
\hom(\mathcal{F}(\Tm_{\class{X}}(1)), \cdot \? ) \colon \class{Y} \to \class{ALG}_{\class{X}}.
\]
In particular, given $\A \in \class{Y}$, the operations on $\hom(\mathcal{F}(\Tm_{\class{X}}(1)), \A )$, for short $\hom(\A)$, are defined as follows:
\[
\psi^{\hom(\A )}(f_{1}, \dots, f_{n}) \coloneqq   \langle f_{1}, \dots, f_{n}\rangle \circ \mathcal{F}(\psi),
\]
for every $f_{1}, \dots, f_{n} \in \hom(\A )$, where $\langle f_{1}, \dots, f_{n}\rangle \colon \mathcal{F}(\Tm_{\class{X}}(n)) \to \A$ is the map induced by the universal property of the coproduct.

 Now, given $\A \in \class{Y}$, we consider the map $\sigma_{\A} \colon \hom(\A) \to \theta_{\mathscr{L}}(\A^{[\kappa]})$ defined by the rule
\[
f \longmapsto \langle f \circ \pi_{1}(x^{i}) : i < \kappa \rangle
\]
where $\pi_{1} \colon \Tm_{\class{Y}}(\kappa) \to \mathcal{F}(\Tm_{\class{X}}(1))$ is the map defined right before Definition \ref{Def:Chap5:ArrowsFunctions}. It is not difficult to see that $\sigma_{\A}$ is a well-defined isomorphism. Hence the global map $\sigma \colon \hom(\mathcal{F}(\Tm_{\class{X}}(1)), \cdot \? ) \to \theta_{\mathscr{L}} \circ [\kappa]$ is a natural isomorphism between $\hom(\mathcal{F}(\Tm_{\class{X}}(1)), \cdot \? ), \theta_{\mathscr{L}} \circ [\kappa] \colon \class{Y} \to \class{ALG}_{\class{X}}$. As a consequence we obtain the following:
\begin{fact}\label{Fact:Chap5:NatIso2}
The map $\hom(\mathcal{F}(\Tm_{\class{X}}(1)), \cdot \? )$ can be viewed as a functor from $\class{Y}$ to $\class{X}$ naturally isomorphic to $\theta_{\mathscr{L}} \circ [\kappa]$.
\end{fact}

Thanks to Facts \ref{Fact:Chap5:natiso1} and \ref{Fact:Chap5:NatIso2}, in order to complete the proof it will be enough to construct a natural isomorphism 
\[
\mu \colon \mathcal{G}^{\ast} \to \hom(\mathcal{F}(\Tm_{\class{X}}(1)), \cdot \? ).
\]
This is what we do now. For every $\A \in \class{Y}$, the component $\mu_{\A}$ of the natural transformation $\mu$ is the following map:
\[
\epsilon_{\A} \circ \mathcal{F}( \cdot ) \? \colon \hom(\Tm_{\class{X}}(1), \mathcal{G}(\A)) \to \hom(\mathcal{F}(\Tm_{\class{X}}(1)), \A).
\]
From the hom-set adjunction associated with $\langle \mathcal{F}, \mathcal{G}, \epsilon, \eta \rangle$ it follows that $\mu_{\A}$ is a bijection. Since $\mathcal{F}$ preserves coproducts, we have that $\mu_{\A}$ is a homomorphism. Therefore we conclude that $\mu_{\A}$ is an isomorphism.

Finally, the fact that the global map $\mu$ satisfies the commutative condition required of natural transformations is witnessed by the hom-set adjunction associated with $\langle \mathcal{F}, \mathcal{G}, \epsilon, \eta \rangle$. Hence $\mu$ is a natural isomorphism as desired.

2. Consider an infinite cardinal $\lambda$ such that $\UUU_{\lambda}(\class{X}) = \class{X}$ and define $\class{K} \coloneqq \GGG\QQQ_{\lambda}(\class{Y}^{[\kappa]})$. Since $\theta$ is compatible with $\mathscr{L}$ in $\class{Y}^{[\kappa]}$ and the compatibility condition is expressible by a set of generalized quasi-equations each of which is written with finitely many variables, we conclude that $\theta$ is compatible with $\mathscr{L}$ in $\class{K}$ too. Moreover, from the fact that $\theta_{\mathscr{L}}(\class{Y}^{[\kappa]}) \subseteq \class{X}$ and $\UUU_{\lambda}(\class{X}) = \class{X}$ it follows that the functor $\theta_{\mathscr{L}} \colon \class{K} \to \class{X}$ is well defined. By Theorems \ref{Thm:McKenzie} and \ref{Thm:ThetaRight} we know that the maps $[\kappa] \colon \class{Y} \to \class{K}$ and $\theta_{\mathscr{L}} 	\colon \class{K} \to \class{X}$ are right adjoint functors. As a consequence their composition $\theta_{\mathscr{L}}\circ [\kappa] \colon \class{Y} \to \class{X}$ is also a right adjoint.
\end{proof}

\begin{Corollary}\label{Cor : Left Adjoint}
Let $\mathcal{F} \colon \class{X} \to \class{Y}$ be a non-trivial left adjoint functor between generalized quasi-varieties and $\phi \in \Con_{\class{X}}\Tm_{\class{X}}(\lambda)$. Assume that the right adjoint of $\mathcal{F}$ decomposes as $\theta_{\mathscr{L}} \circ [\kappa]$. Then
\[
\mathcal{F}(\Tm_{\class{X}}(\lambda) / \phi) \cong \Tm_{\class{Y}}(\kappa\times\lambda) / \textup{Cg}_{\class{Y}}  (\btau^{\ast}(\phi)\cup \bigcup_{j < \lambda}\Theta ( \vec{x}_{j})).
\]
\end{Corollary}

\begin{Remark}\label{Rem:Freyd}
The description of right adjoints given in Theorem \ref{Thm:Main} can be seen as a purely algebraic formulation of the classical description of adjunctions in categories with a free object, which can be traced back at least to \cite{Fr66}.

To see why, suppose that $\mathcal{F} \colon \class{X} \longleftrightarrow \class{Y} \colon \mathcal{G}$ is an adjunction $\mathcal{F} \dashv \mathcal{G}$, and that $\class{X}$ and $\class{Y}$ are prevarieties.\ We proceed to sketch the general description of $\mathcal{G}(\A)$ in \cite{Fr66}. Since $\class{X}$ contains free algebras, the universe of the algebra $\mathcal{G}(\A)$ can be identified with $\hom_{\class{X}}(\Tm_{\class{X}}(1), \mathcal{G}(\A))$.\ By the hom-set adjunction induced by $\mathcal{F} \dashv \mathcal{G}$, we know that
\[
\hom_{\class{X}}(\Tm_{\class{X}}(1), \mathcal{G}(\A)) \cong\hom_{\class{Y}}(\mathcal{F}(\Tm_{\class{X}}(1)), \A).
\]
Since $\class{Y}$ contains arbitrarily large free algebras, the algebra $\mathcal{F}(\Tm_{\class{X}}(1))$ can be expressed as a suitable quotient of a free algebra, i.e.\ $\mathcal{F}(\Tm_{\class{X}}(1)) \cong \Tm_{\class{Y}}(\kappa)/ \theta$ for some cardinal $\kappa$ and some congruence $\theta$. Thus the universe of $\mathcal{G}(\A)$ can be identified with $\hom_{\class{Y}}(\Tm_{\class{Y}}(\kappa)/ \theta, \A)$. More in general $\mathcal{G}(\A)$ can be identified with the set $\hom_{\class{Y}}(\Tm_{\class{Y}}(\kappa)/ \theta, \A)$, equipped with a suitable algebraic structure. This provides a full arrow-theoretic description of the algebra $\mathcal{G}(\A)$ as
\[
\mathcal{G}(\A) \cong \hom_{\class{Y}}(\Tm_{\class{Y}}(\kappa)/ \theta, \A).
\]

The main contribution of the present work is to recognize that the algebra $\hom_{\class{Y}}(\Tm_{\class{Y}}(\kappa)/ \theta, \A)$ in the above display can be given a very transparent description in terms of matrix powers and compatible equations.
%
% This is a consequence of the fact that the set $\hom_{\class{Y}}(\Tm_{\class{Y}}(\kappa)/ \theta, \A)$ can be identified with the set of solutions of the equations $\theta(\vec{x})$ in $\A$, which is exactly the universe of the algebra $\theta_{\mathscr{L}}(\A^{[\kappa]})$. Moreover, this identification respect the algebraic structures, yielding an isomorphism
%\[
%\hom_{\class{Y}}(\Tm_{\class{Y}}(\kappa)/ \theta, \A) \cong \theta_{\mathscr{L}}(\A^{[\kappa]}).
%\]
%As a consequence the structure of the algebra $\mathcal{G}(\A)$ can be expressed in purely algebraic and combinatorial terms as $\theta_{\mathscr{L}}(\A^{[\kappa]})$. In particular, this description of $\mathcal{G}(\A)$ was exploited to establish the correspondence between between adjunctions and contextual translations.
\qed
\end{Remark}

Until now we showed that every right adjoint functor $\mathcal{G} \colon \class{Y} \to \class{X}$ between generalized quasi-varieties induces a contextual translation $\langle \btau, \theta\rangle$ of $\vDash_{\class{X}}$ into $\vDash_{\class{Y}}$, and vice-versa. In general, contextual translations $\langle \btau, \theta\rangle$ are infinite objects, in the sense that $\btau$ is a map that translates terms into possibly infinite sequences of terms and $\theta$ is a possibly infinite set of equations. It is therefore natural to ask under which conditions these contextual translations can be finitized. The next lemma provides an answer in the case where $\class{X}$ and $\class{Y}$ are quasi-varieties.

\begin{Lemma}\label{Lem : Fin Pres}
Let $\mathcal{F} \colon \class{X} \longleftrightarrow \class{Y} \colon \mathcal{G}$ be an adjunction $\mathcal{F} \dashv \mathcal{G}$ between quasi-varieties. The following conditions are equivalent:
\benroman
\item $\mathcal{F}$ preserves finitely presentable algebras.
\item $\mathcal{F}(\Tm_{\class{X}}(1))$ is finitely presentable.
\item $\mathcal{G}$ preserves directed colimits.
\item $\mathcal{G}$ can be decomposed as $\theta_{\mathscr{L}} \circ [\kappa]$  with both $\kappa$ and $\theta$ finite.
\eroman
\end{Lemma}
\begin{proof}
The equivalence between (i) and (iii) is well known, and is a consequence of the fact that the finitely $\class{X}$-presentable algebras are exactly the algebras $\A \in \class{X}$ for which the functor $\hom(\A,\cdot) \colon\A \to \class{Set}$ preserves directed colimits (see Lemma \ref{Lem : What are presetable}). Part (i)$\Rightarrow$(ii) is trivial and part (iv)$\Rightarrow$(i) is a consequence of Corollary \ref{Cor : Left Adjoint}.

(ii)$\Rightarrow$(iv): Assume that $\mathcal{F}(\Tm_{\class{X}}(1))$ is finitely presentable. Then there are $n \in \omega$ and a compact $\class{Y}$-congruence $\Theta$ such that $\mathcal{F}(\Tm_{\class{X}}(1)) = \Tm_{\class{Y}}(n)/ \Theta$. Now, $\Theta$ is generated by a finite set $\Phi \subseteq \Theta$. This means that 
$\mathcal{G}$ can be decomposed as $\theta_{\mathscr{L}} \circ [n]$, where $\theta \coloneqq \{ \vec{\epsilon} \thickapprox \vec{\delta} : \langle \epsilon, \delta\rangle \in \Phi \}$ and $\vec{\epsilon}$, $\vec{\delta}$ are sequences of length $n$.
\end{proof}

%\begin{Remark}
%Even if we do not pursue the details here, it is interesting to observe that if $\mathcal{G} \colon \class{Y} \to \class{X}$ is a right adjoint between quasi-varieties satisfying any of the equivalent conditions in the above lemma, then $\mathcal{G}$ is indeed induced by a \textit{model theoretic} interpretation of the language of $\class{X}$ into the language of $\class{Y}$, see \cite[Chapter 5]{Ho93} for the relevant definitions.
%\qed
%\end{Remark}

The next example shows that there are adjunctions between quasi-varieties that do not meet the equivalent conditions of Lemma \ref{Lem : Fin Pres}. In other words, it shows that there are contextual translations between finitary relative equational consequences that cannot be finitized.

\begin{exa}[\textsf{Ring Hom-Functor}]\label{Ex:Chap6:Hom-Functor}
Consider a generalized quasi-variety  $\class{X}$ and an algebra $\A \in \class{X}$. Then let $\hom(\A, \cdot) \colon \class{X} \to \class{Set}$ be the functor defined by the following rule:
\begin{align*}
\B &\longmapsto \hom(\A, \B)\\
f \colon \B \to \C &\longmapsto f \circ (\cdot) \colon \hom(\A, \B) \to \hom(\A, \C).
\end{align*}
The functor $\hom(\A, \cdot)$ has a left adjoint $\mathcal{F} \colon \class{Set} \to \class{X}$ defined as follows. Given a set $I$, the algebra $\mathcal{F}(I)$ is the copower of $\A$ indexed by $I$. Moreover, given a function $f \colon I \to J$ between sets, we let $\mathcal{F}(f)\colon\mathcal{F}(I) \to \mathcal{F}(J)$ be the map $\langle p_{f(i)} : i \in I \rangle$ induced by the universal propery of the coproduct $\mathcal{F}(I)$, where $\{ p_{j} \colon \A \to \mathcal{F}(J) : j \in J \}$ are the maps associated with the copower $\mathcal{F}(J)$.

Now consider the special case where $\class{X}$ is the variety $\class{R}$ of commutative rings with unit. Then consider the functor $\mathcal{F}$ that is left adjoint to $\hom(\boldsymbol{Q}, \cdot) \colon \class{R} \to \class{Set}$, where $\boldsymbol{Q}$ is the ring of rational numbers. First observe that $\mathcal{F}$ does not preserve finitely generated algebras. Observe that finitely generated algebras are exactly the quotients of the finitely presentable ones. Since $\mathcal{F}$ preserves surjective homomorphisms, we conclude that it does not preserve finitely presentable algebras.
\qed
\end{exa}

%In this paper contextual translations have been presented as translations between \textit{relative equational consequences}. Nevertheless, the contextual translations coming from some of the motivating examples (such as G\"odel and Kolmogorov's ones) originated as translations between \textit{propositional logics}. More precisely, it is a general fact that if two generalized quasi-varieties $\class{X}$ and $\class{Y}$ are the equivalent algebraic semantics of two algebraizable logics $\LL$ and $\LL'$ in the sense of \cite{BP89}, then every contextual translation of $\vDash_{\class{X}}$ into $\vDash_{\class{Y}}$ can be viewed as translation of $\LL$ into $\LL'$. In this sense, contextual translation may provide a useful notion of translation between algebraizable logics (possibly in different languages).

\paragraph{\bfseries Acknowledgements.}
I am grateful to James Raftery for raising the question of whether category equivalence may provide a notion of equivalence between propositional logics. Thanks are due also to Josep Maria Font, Ramon Jansana and Juan Climent, who read carefully many versions of this work and provided several useful comments that improved the readability of this paper.  I wish to thank also Ji\v{r}\'i Velebil and Mat\v{e}j Dost\'al for some useful observations on categorical universal algebra. Finally, thanks are due to the anonymous referee for providing additional references and useful remarks, which helped to improve the presentation.  This research was supported by the joint project of Austrian Science Fund (FWF) I$1897$-N$25$ and Czech Science Foundation (GACR) $15$-$34650$L

\end{document}